\documentclass[11pt]{article}

\usepackage{amsfonts}
\usepackage{amssymb}
\usepackage{amsthm}
\usepackage{amstext}
\usepackage{amscd}
\usepackage{amsmath}
\usepackage{delarray}

\usepackage{epsfig}

\usepackage{fullpage}

\usepackage{color}
\definecolor{light-blue}{rgb}{0.8,0.85,1}
\definecolor{blue}{rgb}{0,0,1}
\definecolor{red}{rgb}{1,0,0}

\DeclareMathOperator{\Res}{Res}

\begin{document}

\title{\bf Phase transition of the largest eigenvalue for
non-null complex sample covariance matrices}

\author{{\bf Jinho Baik}
\footnote{Department of Mathematics, University of Michigan, Ann
Arbor, MI, 48109, USA, baik@umich.edu}, {\bf G\'erard Ben
Arous}\footnote{Department of Mathematics, Courant Institute of
Mathematical Sciences, New York, NY, 10012, USA,
benarous@cims.nyu.edu} and {\bf Sandrine
P\'ech\'e}\footnote{Department of Mathematics Ecole Polyechnique
F\'ed\'erale de Lausanne, 1015 Lausanne Switzerland,
sandrine.peche@epfl.ch; Current address:  Institut Fourier, UJF
Grenoble 38000 France,
sandrine.peche@ujf-grenoble.fr}}

\date{\today}

\maketitle
\newtheorem{theo}{Theorem}[section]
\newtheorem{prop}{Proposition}[section]
\newtheorem{lemma}{Lemma}[section]
\newtheorem{cor}{Corollary}[section]
\newtheorem{conjecture}{Conjecture}[section]
\newtheorem{definition}{Definition}[section]
\newcommand{\bconj}{\begin{conjecture}}
\newcommand{\econj}{\end{conjecture}}
\newcommand{\bdefi}{\begin{definition}}
\newcommand{\edefi}{\end{definition}}
\newcommand{\bt}{\begin{theo}}
\newcommand{\Si}{\Sigma}
\newcommand{\et}{\end{theo}}
\newcommand{\bp}{\begin{prop}}
\newcommand{\ep}{\end{prop}}
\newcommand{\bl}{\begin{lemme}}
\newcommand{\el}{\end{lemme}}
\newcommand{\be}{\begin{equation}}
\newcommand{\ee}{\end{equation}}
\newcolumntype{L}{>{$}l<{$}}
\newenvironment{Cases}{\begin{array}\{{lL.}}{\end{array}}

\newcommand{\qq}{q}

\theoremstyle{remark}
\newtheorem*{rem}{Remark}
\newtheorem{remark}{Remark}[section]
\newtheorem*{nrem}{\textbf{Notational Remark}}
\newtheorem*{Conjecture}{\textbf{Conjecture}}

\begin{abstract}
We compute the limiting distributions of the largest eigenvalue of
a complex Gaussian sample covariance matrix when both the number
of samples and the number of variables in each sample become
large. When all but finitely many, say $r$, eigenvalues of the
covariance matrix are the same, the dependence of the limiting
distribution of the largest eigenvalue of the sample covariance
matrix on those distinguished $r$ eigenvalues of the covariance
matrix is completely characterized in terms of an infinite
sequence of new distribution functions that generalize the
Tracy-Widom distributions of the random matrix theory. Especially
a phase transition phenomena is observed. Our results also apply
to a last passage percolation model and a queuing model.
\end{abstract}

\section{Introduction}

Consider $M$  independent, identically distributed samples
$\vec{y}_1,\dots, \vec{y}_M$, all of which are $N\times 1$ column
vectors. We further assume that the sample vectors $\vec{y}_k$ are
Gaussian with mean $\vec{\mu}$ and \emph{covariance} $\Si$, where
$\Si$ is a fixed $N\times N$ positive matrix; the density of a
sample $\vec{y}$ is
\begin{equation}\label{eq:sampledensity}
  p(\vec{y})= \frac{1}{(2\pi)^{N/2} (\det \mathbf{\Sigma})^{1/2}}
e^{-\frac12 <\vec{y}-\vec{\mu},
\mathbf{\Sigma}^{-1}(\vec{y}-\vec{\mu})>} ,
\end{equation}
where $<,>$ denotes the inner product of vectors. We denote by
$\ell_1,\dots, \ell_N$ the eigenvalues of the covariance matrix
$\Sigma$, called the `population eigenvalues'. The sample mean
$\overline{Y}$ is defined by
$\overline{Y}:=\frac1{M}(\vec{y}_1+\cdots+\vec{y}_M)$ and we set
$X=[\vec{y}_1-\overline{Y},\cdots, \vec{y}_M-\overline{Y}]$ to be
the (centered) $N\times M$ sample matrix. Let $S=\frac1{M} XX'$ be
the \emph{sample covariance matrix}. The eigenvalues of $S$,
called the `sample eigenvalues', are denoted by $\lambda_1
>\lambda_2
> \cdots
>\lambda_N>0$. (The eigenvalues are simple with probability $1$.)
The probability space of $\lambda_j$'s is sometimes called the
Wishart ensemble (see e.g. \cite{Muirhead}).

Contrary to the traditional assumptions, it is of current interest
to study the case when $N$ is of same order as $M$. Indeed when
$\Sigma=I$ (null-case), several results are known. As
$N,M\to\infty$ such that $M/N\to \gamma^2 \ge 1$, the following
holds.

\begin{itemize}
\item \textbf{Density of eigenvalues} \cite{MP}: For any real $x$,
\begin{equation}\label{eq:DOSreal}
  \frac1{N} \#\{ \lambda_j : \lambda_j\le x\} \to H(x)
\end{equation}
where
\begin{equation}\label{eq:DOSformula}
  H'(x) = \frac{\gamma^2}{2\pi x}\sqrt{(b-x)(x-a)}, \quad a<x<b,
\end{equation}
and $a=(\frac{\gamma-1}{\gamma})^2$ and
$b=(\frac{1+\gamma}{\gamma})^2$.
\item \textbf{Limit of the largest eigenvalue} \cite{Geman}:
\begin{equation}\label{eq:SLreal}
  \lambda_1 \to \biggl(\frac{1+\gamma}{\gamma}\biggr)^2 \quad \text{a.s.}
\end{equation}
\item \textbf{Limiting distribution}
\cite{Johnstone}: For any real $x$,
\begin{equation}\label{eq:CLTreal}
  \mathbb{P} \biggl( \bigr( \lambda_1 - \bigl(\frac{\gamma+1}\gamma \bigr)^2
  \bigr) \cdot
  \frac{\gamma M^{2/3}}{(1+\gamma)^{4/3}} \le x \biggr)
  \to F_{GOE}(x)
\end{equation}
where $F_{GOE}(x)$ is the so-called GOE Tracy-Widom distribution,
which is the limiting distribution of the largest eigenvalue of a
random real symmetric matrix from the Gaussian orthogonal ensemble
(GOE) as the size of the matrix tends to infinity
\cite{TracyWidomGOE}.
\item \textbf{Robustness to models} \cite{Soshnikovcovariance}: It
turned out that the Gaussian assumption is unnecessary and a
result similar to \eqref{eq:CLTreal} still holds for a quite
general class of independent, identically distributed random
samples.
\end{itemize}

From \eqref{eq:DOSreal} and \eqref{eq:SLreal}/\eqref{eq:CLTreal},
we find that the largest sample eigenvalue $\lambda_1$ in the null
case converges to the rightmost edge of support of the limiting
density of eigenvalues. However, in practice (see e.g.
\cite{Johnstone}) there often are statistical data for which one
or several large sample eigenvalues are separated from the bulk of
the eigenvalues. For instance, see Figure 1 and 2 of the paper
\cite{Johnstone} which plot the sample eigenvalues of the
functional data consisting of a speech dataset of 162 instances of
a phoneme ``dcl'' spoken by males calculated at 256 points
\cite{Buja}. Other examples of similar phenomena include
mathematical finance \cite{PlerousGRAGS}, \cite{LalouxCPB},
\cite{MalevergneS}, wireless communication \cite{Telatar}, physics
of mixture \cite{SearC}, and data analysis and statistical
learning \cite{HoyleR}. As suggested in \cite{Johnstone}, this
situation poses a natural question: when $\Sigma\neq I$ (non-null
case), how do a few large sample eigenvalues depend on the
population eigenvalues? More concretely, if there are a few large
population eigenvalues, do they pull to the sample eigenvalues,
and for it to happen, how large the population eigenvalues should
be?

Though this might be a challenging problem for real sample data,
it turned out that one could answer some of the above questions in
great detail for complex Gaussian samples. Complex sample
covariance matrix has an application in multi-antenna Gaussian
channels in wireless communication \cite{Telatar}. Also the
results of complex case lead us to a guess for aspects of the real
case (see Conjecture in section \ref{sec:mainresults} below).
Another reason of studying complex sample covariance matrix is its
relation to a last passage percolation model and a queueing
theory. See section \ref{sec:last} below for such a connection.

Before we present our work, we first summarize some known results
for the complex sample covariance matrices.

\subsection{Some known results for the eigenvalues of complex sample
covariance matrices}

We assume that the samples $\vec{y}$ are \emph{complex Gaussian}
with mean $\vec{\mu}$ and covariance $\Sigma$. Hence the density
of $\vec{y}$ is precisely given by \eqref{eq:sampledensity} with
the understanding that $<,>$ denotes now the complex inner
product. The (centered) sample matrix $X$ and the sample
covariance matrix $S=\frac1{N}XX^*$ are defined as before where
$X^*$ is the transpose followed by the complex conjugation. Recall
that the eigenvalues of $S$, sample eigenvalues, are denoted by
$\lambda_1\ge \cdots \ge \lambda_N>0$, and the eigenvalues of
$\Sigma$, population eigenvalues, are denoted by $\ell_1, \dots,
\ell_N>0$.

\begin{itemize}
\item \textbf{Density of eigenvalues} \cite{MP}, \cite{BaiSilverstein}
(see also Theorem 3.4 of \cite{BaiReview}): When all but finitely
many eigenvalues $\ell_j$ of $\Sigma$ are equal to $1$, as
$M,N\to\infty$ such that $M/N\to \gamma^2\ge 1$, the limiting
density of the sample eigenvalues $\lambda_j$ is given by
\begin{equation}\label{eq:DOS}
  \frac1{N} \#\{ \lambda_j : \lambda_j\le x\} \to H(x)
\end{equation}
where $H(x)$ is again defined by \eqref{eq:DOSformula}.
\item \textbf{Null case} : When $\Sigma=I$, as $M,N\to\infty$
such that $M/N\to \gamma^2\ge 1$, \cite{Geman}
\begin{equation}\label{eq:SL}
  \lambda_1 \to \biggl(\frac{1+\gamma}{\gamma}\biggr)^2 \quad \text{a.s.}
\end{equation}
and for any real $x$ (see, e.g. \cite{Forrester},
\cite{JohanssonOrthogonal})
\begin{equation}\label{eq:CLTnull}
  \mathbb{P} \biggl( \bigr( \lambda_1 - \bigl(\frac{1+\gamma}{\gamma}\bigr)^2 \bigr) \cdot
  \frac{\gamma M^{2/3}}{(1+\gamma)^{4/3}} \le x \biggr)
  \to F_{GUE}(x)
\end{equation}
where $F_{GUE}(x)$ is the GUE Tracy-Widom distribution, which is
the limiting distribution of the largest eigenvalue of a random
complex Hermitian matrix from the Gaussian unitary ensemble (GUE)
as the size of the matrix tends to infinity \cite{TracyWidom}.
Moreover, the limit \eqref{eq:CLTnull} holds true for a quite
general class of independent, identically distributed random
samples, after suitable scaling \cite{Soshnikovcovariance}.
\end{itemize}

\begin{remark}
The distribution function $F_{GUE}$ is different from $F_{GOE}$. A
formula of $F_{GUE}(x)$ is given in \eqref{eq:TWGUE} below and a
formula for $(F_{GOE}(x))^2$ is given in \eqref{eq:GOEPII} below.
\end{remark}

\begin{remark}
When $\Sigma=I$, the probability space of the eigenvalues
$\lambda_j$ of $S$ is sometimes called the Laguerre unitary
ensemble (LUE) since the correlation functions of $\lambda_j$ can
be represented in terms of Laguerre polynomials. Similarly, for
real samples with $\Sigma=I$, the probability space of the
eigenvalues of $S$ is called the Laguerre orthogonal ensemble
(LOE). See e.g. \cite{ForresterBook}.
\end{remark}

Note that the limiting density of the eigenvalues $\lambda_j$ is known for
general
$\Sigma \neq I$, but
the convergence \eqref{eq:SL}/\eqref{eq:CLTnull} of $\lambda_1$ to the edge of
the support of the limiting distribution of the eigenvalues
was obtained only when $\Sigma=I$.
The following result of P{\'e}ch{\'e} \cite{Sandrine} generalizes
\eqref{eq:CLTnull}
and shows that when
all but finitely many eigenvalues $\ell_k$ of $\Sigma$ are $1$
and those distinguished eigenvalues are ``not too big'',
$\lambda_1$ is still not separated from the rest of the eigenvalues.

\begin{itemize}
\item
When $\ell_{r+1}=\dots=\ell_N=1$ for a fixed integer
$r$ and $\ell_1=\dots=\ell_r<2$ are fixed,
 as $M=N\to\infty$,
\cite{Sandrine}
\begin{equation}
  \mathbb{P} \bigl( \bigr( \lambda_1 - 4 \bigr) \cdot
  2^{-4/3}M^{2/3} \le x \bigr)
  \to F_{GUE}(x).
\end{equation}
\end{itemize}

A natural question is then whether the upper bound 2 of
$\ell_1=\cdots =\ell_r$ is critical. One of our result in this
paper is that it is indeed the critical value. Moreover, we find
that if some of $\ell_j$ are precisely equal to the critical
value, then the limiting distribution is changed to something new.
And if one or more $\ell_j$ are bigger than the critical value,
the fluctuation order $M^{2/3}$ is changed to the Gaussian type
order $\sqrt{M}$. In order to state our results, we first need
some definitions.


\subsection{Definitions of some distribution functions}

\subsubsection{Airy-type distributions}

Let $Ai(u)$ be the Airy function which has the integral
representation
\begin{equation}\label{eq:Airyintrep}
  Ai(u) = \frac{1}{2\pi } \int e^{iua+i\frac{1}3 a^{3}}da
\end{equation}
where the contour is from $\infty e^{5i\pi/6}$ to $\infty
e^{i\pi/6}$. Define the \emph{Airy kernel} (see e.g.
\cite{TracyWidom}) by
\begin{equation}
  \mathbf{A}(u,v) = \frac{Ai(u)Ai'(v)-Ai'(u)Ai(v)}{u-v}
\end{equation}
and let $\mathbf{A}_x$ be the operator acting on $L^2((x,\infty))$
with kernel $\mathbf{A}(u,v)$. An alternative formula of the Airy
kernel is
\begin{equation}\label{eq:Airykernelint}
  \mathbf{A}(u,v) = \int_0^\infty Ai(u+z)Ai(z+v) dz,
\end{equation}
which can be checked directly by using the relation
$Ai''(u)=uAi(u)$ and integrating by parts. For $m=1,2,3,\dots$,
set
\begin{equation}\label{eq:sdef}
  s^{(m)}(u) = \frac{1}{2\pi }\int e^{iua+ i\frac13a^3} \frac1{(ia)^m} da
\end{equation}
where the contour is from $\infty e^{5i\pi/6}$ to $\infty
e^{i\pi/6}$ such that the point $a=0$ lies \emph{above} the
contour. Also set
\begin{equation}
  t^{(m)}(v) = \frac{1}{2\pi } \int e^{iva+ i\frac13 a^3} (-ia)^{m-1}
  da
\end{equation}
where the contour is from $\infty e^{5i\pi/6}$ to $\infty
e^{i\pi/6}$. Alternatively,
\begin{equation}
  s^{(m)}(u) =  \sum_{\substack{\ell+3n = m-1\\ \ell, n = 0,1,2,\dots}} \biggl\{ \frac{(-1)^n}{3^n \ell ! n!}
  u^\ell + \frac1{(m-1)!} \int_\infty^u (u-y)^{m-1} Ai(y)dy
  \biggr\}
\end{equation}
and
\begin{equation}
  \qquad t^{(m)}(v) = \biggl( -\frac{d}{dv} \biggr)^{m-1} Ai(v).
\end{equation}
See Lemma \ref{lem:Fk} below for the proof that the two formulas of $s^{(m)}(u)$
are the same.

\begin{definition}\label{def:Fk}
For $k=1,2,\dots$, define for real $x$,
\begin{equation}\label{eq:Fkdef}
  F_k(x) = \det (1-\mathbf{A}_x) \cdot
  \det \biggl( \delta_{mn} - < \frac{1}{1-\mathbf{A}_x}
  s^{(m)}, t^{(n)} > \biggr)_{1\le m,n\le k},
\end{equation}
where $<,>$ denotes the (real) inner product of functions in
$L^2((x,\infty))$. Let $F_0(x)=\det(1-\mathbf{A}_x)$.
\end{definition}

The fact that the inner product in \eqref{eq:Fkdef} makes sense
and hence $F_k(x)$ is well-defined is proved in Lemma \ref{lem:Fk}
below.

It is well-known that (see e.g. \cite{Forrester},
\cite{TracyWidom})
\begin{equation}\label{eq:TWGUE}
  F_0(x)= \det(1-\mathbf{A}_x)=F_{GUE}(x)
\end{equation}
and hence $F_0$ is the GUE Tracy-Widom distribution function.
There is an alternative expression of $F_0$. Let $u(x)$ be the
solution to the Painlev\'e II equation
\begin{equation}
  u''=2u^3+xu
\end{equation}
satisfying the condition
\begin{equation}
  u(x) \sim -Ai(x), \qquad x\to +\infty.
\end{equation}
There is a unique, global solution \cite{HM}, and satisfies (see
e.g. \cite{HM}, \cite{DZ2})
\begin{eqnarray}
  u(x) &=& -\frac{e^{-\frac23x^{3/2}}}{2\sqrt{\pi}x^{1/4}}
  + O\biggl( \frac{e^{-\frac43x^{3/2}}}{x^{1/4}} \biggr) \qquad \text{as
  $x\to+\infty$} \\
  u(x) &=& -\sqrt{\frac{-x}{2}} \bigl( 1+ O(x^{-2}) \bigr)
  \qquad \text{as $x\to -\infty$}.
\end{eqnarray}
Then \cite{TracyWidom}
\begin{equation}\label{eq:FredPII}
  F_0(x)= \det(1-\mathbf{A}^{(0)}_x) = \exp \biggl( -\int_x^\infty (y-x)
  u^2(y)dy \biggr).
\end{equation}
In addition to being a beautiful identity, the right-hand-side of
\eqref{eq:FredPII} provides a practical formula to plot the graph
of $F_0$.

For $k=1$, it is known that (see \cite{ForresterSamll}, (3.34) of
\cite{BorodinForrester})
\begin{equation}\label{eq:GOEPII}
  F_1(x)= \det(1-\mathbf{A}_x)\cdot
  \biggr( 1- <\frac1{1-\mathbf{A}_x} s^{(1)}, t^{(1)} >\biggr)
  = (F_{GOE}(x))^2 .
\end{equation}
The function $F_{GOE}$ also has a Painlev\'e formula
\cite{TracyWidomGOE} and
\begin{equation}\label{eq:F1PII}
  F_1(x) = F_0(x) \exp\bigg(\int_x^\infty u(y)dy \biggr).
\end{equation}

The functions $F_k$, $k\ge 2$, seem to be new. The Painlev\'e
formula of $F_k$ for general $k\ge 2$ will be presented in
\cite{Baik}. For each $k\ge 2$, $F_k(x)$ is clearly a continuous
function in $x$. Being a limit of non-decreasing functions as
Theorem \ref{thm:main} below shows, $F_k(x)$ is a non-decreasing
function. It is also not difficult to check by using a
steepest-descent analysis that $F_k(x) \to 1$ as $x\to +\infty$
(cf. Proof of Lemma \ref{lem:Fk}). However, the proof that
$F_k(x)\to 0$ as $x\to -\infty$ is not trivial. The fact that
$F_k(x)\to 0$ as $x\to -\infty$ is obtained in \cite{Baik} using
the Painlev\'e formula. Therefore $F_k(x)$, $k\ge 2$, are
distribution functions, which generalize the Tracy-Widom
distribution functions. (The functions $F_0, F_1$ are known to be
distribution functions.)

\subsubsection{Finite GUE distributions}

Consider the density of $k$ particles $\xi_1,\dots, \xi_k$ on the
real line defined by
\begin{equation}
  p(\xi_1,\dots,\xi_k)
  = \frac1{Z_k} \prod_{1\le i<j\le k} |\xi_i-\xi_j|^2
  \cdot \prod_{i=1}^k e^{-\frac12 \xi_i^2}
\end{equation}
where $Z_k$ is the normalization constant,
\begin{equation}
  Z_k := \int_{-\infty}^\infty \cdots \int_{-\infty}^\infty
  \prod_{1\le i<j\le k} |\xi_i-\xi_j|^2
  \cdot \prod_{i=1}^k e^{-\frac12 \xi_i^2}
  d\xi_1\cdots d\xi_k
  = (2\pi)^{k/2} \prod_{j=1}^k j!,
\end{equation}
which is called the Selberg's integral (see e.g. \cite{Mehta}).
This is the density of the eigenvalues of the Gaussian unitary
ensemble (GUE), the probability space of $k\times k$ Hermitian
matrices $H$ whose entries are independent Gaussian random
variables with mean $0$ and standard deviation $1$ for the
diagonal entries, and mean $0$ and standard deviation $1/2$ for
each of the real and complex parts of the off-diagonal entries
(see e.g. \cite{Mehta}).

\begin{definition}
For $k=1,2,3,\dots$, define the distribution $G_k(x)$ by
\begin{equation}\label{eq:Gkdef}
  G_k(x) = \frac1{Z_k}  \int_{-\infty}^x\cdots \int_{-\infty}^x
  \prod_{1\le i<j\le k} |\xi_i-\xi_j|^2
  \cdot \prod_{i=1}^k e^{-\frac12 \xi_i^2}
  d\xi_1\cdots d\xi_k .
\end{equation}
\end{definition}

In other words, $G_k$ is the distribution of the \emph{largest
eigenvalue} of $k\times k$ GUE. When $k=1$, this is the Gaussian
distribution,
\begin{equation}
  G_1(x) = \frac1{\sqrt{2\pi}} \int_{-\infty}^x
  e^{-\frac12 \xi_1^2} d\xi_1
  = erf(x) .
\end{equation}
There is an alternative expression of $G_k$ in terms of a Fredholm
determinant similar to the formula \eqref{eq:Fkdef} of $F_k$. Let
$p_n(x)=c_nx^n+\cdots$ be the polynomial of degree $n$ ($c_n> 0$
is the leading coefficient) determined by the orthogonality
condition
\begin{equation}
  \int_{-\infty}^\infty p_m(x)p_n(x) e^{-\frac12 x^2} dx =
  \delta_{mn}.
\end{equation}
The orthonormal polynomial $p_n$ is given by
\begin{equation}\label{eq:porthonomal}
  p_n(\xi):= \frac1{(2\pi)^{1/4}2^{n/2}\sqrt{n!}}
  H_n(\frac{\xi}{\sqrt{2}}),
\end{equation}
where $H_n$ is the Hermite polynomial. The leading coefficient
$c_n$ of $p_n$ is (see e.g. \cite{Koekoek})
\begin{equation}
  c_n= \frac{1}{(2\pi)^{1/4}\sqrt{n!}} .
\end{equation}
Then the so-called orthogonal polynomial method in the random
matrix theory establishes that:
\begin{lemma}
For any $k=1,2,\dots$ and $x\in\mathbb{R}$,
\begin{equation}\label{eq:Hop}
  G_k(x) = \det(1-\mathbf{H}_x^{(k)}),
\end{equation}
where $\mathbf{H}^{(k)}_x$ is the operator acting on
$L^2((x,\infty))$ defined by the kernel
\begin{equation}
  \mathbf{H}^{(k)}(u,v) = \frac{c_{k-1}}{c_{k}}
  \frac{p_k(u)p_{k-1}(v)-p_{k-1}(u)p_k(v)}{u-v} e^{-(u^2+v^2)/4}.
\end{equation}
\end{lemma}

This is a standard result in the theory of random matrices. The
proof can be found, for example, in \cite{Mehta}, \cite{TW2}.
There is also an identity of the form \eqref{eq:FredPII} for
$\det(1-\mathbf{H}_x^{(k)})$, now in terms of Painlev\'e IV
equation. See \cite{TracyWidomFred}.

\subsection{Main Results}\label{sec:mainresults}

We are now ready to state our main results.

\begin{theo}\label{thm:main}
Let $\lambda_1$ be the largest eigenvalue of the sample covariance
matrix constructed from $M$ independent, identically distributed
complex Gaussian sample vectors of $N$ variables. Let $\ell_1,
\cdots, \ell_N$ denote the eigenvalues of the covariance matrix of
the samples. Suppose that for a fixed integer $r\ge 0$,
\begin{equation}
  \ell_{r+1}=\ell_{r+2}=\dots =\ell_N=1.
\end{equation}
As $M,N\to \infty$ while $M/N=\gamma^2$ is in a compact subset of
$[1,\infty)$, the following holds for any real $x$ in a compact
set.
\begin{itemize}
\item[(a)] When for some $0\le k\le r$,
\begin{equation}
  \ell_1=\dots= \ell_k=1+\gamma^{-1}
\end{equation}
and $\ell_{k+1}, \dots, \ell_{r}$ are in a compact subset
of $(0, 1+\gamma^{-1})$,
\begin{equation}\label{eq:mainF}
  \mathbb{P} \biggl( \bigr( \lambda_1 - (1+\gamma^{-1})^2 \bigr) \cdot
  \frac{\gamma}{(1+\gamma)^{4/3}} M^{2/3} \le x \biggr)
  \to F_k(x).
\end{equation}
where $F_k(x)$ is defined in \eqref{eq:Fkdef}.
\item[(b)] When for some $1\le k\le r$,
\begin{equation}
  \text{$\ell_1=\cdots=\ell_k$ is in a compact set of $(1+\gamma^{-1}, \infty)$}
\end{equation}
and $\ell_{k+1}, \dots, \ell_r$ are in a compact subset of
$(0, \ell_1)$,
\begin{equation}\label{eq:mainG}
  \mathbb{P} \biggl( \bigr( \lambda_1 -
  \bigl(\ell_1+\frac{\ell_1\gamma^{-2}}{\ell_1-1}\bigr) \bigr) \cdot
  \frac{\sqrt{M}}{\sqrt{\ell_1^2-\frac{\ell_1^2\gamma^{-2}}{(\ell_1-1)^2}}} \le x \biggr)
  \to G_k(x)
\end{equation}
where $G_k(x)$ is defined in \eqref{eq:Gkdef}.
\end{itemize}
\end{theo}

\begin{figure}[ht]
\centerline{\epsfxsize=8cm\epsfbox{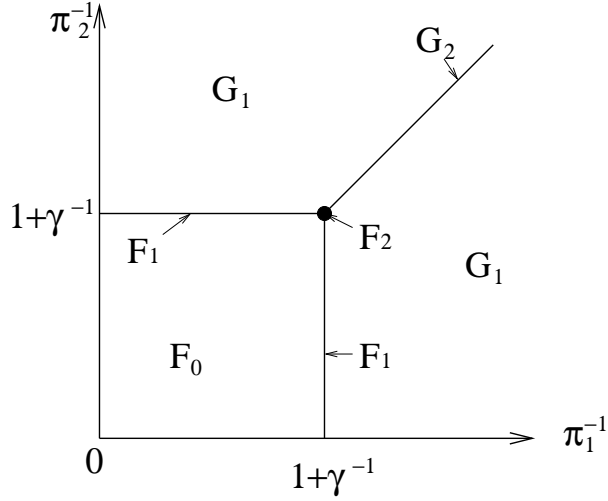}}
\caption{Diagram of the limiting distributions for various choices
of $\ell_1=\pi_1^{-1}$ and $\ell_2=\pi_2^{-1}$ while
$\ell_3=\cdots=\ell_N=1$.} \label{fig:Phasediag}
\end{figure}

Hence, for instance, when $r=2$,
\begin{equation}
  \ell_3=\cdots =\ell_N=1
\end{equation}
and there are two distinguished eigenvalues $\ell_1$ and $\ell_2$
of the covariance matrix. Assume without loss of generality that
$\ell_1\ge \ell_2$. Then
\begin{equation}
  \mathbb{P} \biggl( \bigr( \lambda_1 - (1+\gamma^{-1})^2 \bigr) \cdot
  \frac{\gamma}{(1+\gamma)^{4/3}} M^{2/3} \le x \biggr)
  \to \begin{cases}  F_0(x), \qquad &
  0< \ell_1, \ell_2 <1+\gamma^{-1} \\
  F_1(x), \qquad & 0<\ell_2 <1+\gamma^{-1}=\ell_1 \\
  F_2(x), \qquad & \ell_1=\ell_2=1+\gamma^{-1} \\
  \end{cases}
\end{equation}
and
\begin{equation}
\begin{split}
  & \mathbb{P} \biggl( \bigr( \lambda_1 -
  \bigl(\ell_1+\frac{\ell_1\gamma^{-2}}{\ell_1-1}\bigr) \bigr) \cdot
  \frac{\sqrt{M}}{\sqrt{\ell_1^2-\frac{\ell_1^2\gamma^{-2}}{(\ell_1-1)^2}}} \le x
  \biggr)
  \to \begin{cases}
  G_1(x), \qquad & \ell_1 > 1+\gamma^{-1}, \, \ell_1>\ell_2  \\
  G_2(x), \qquad & \ell_1=\ell_2> 1+\gamma^{-1},
  \end{cases}
\end{split}
\end{equation}
assuming that $\ell_1, \ell_2$ are in compact sets in each case.
See Figure \ref{fig:Phasediag} for a diagram.

Note the different fluctuation orders $M^{2/3}$ and $\sqrt{M}$
depending on the values of $\ell_1, \ell_2$. This type of `phase
transition' was also observed in \cite{BR2, BR4, SpohnPCurrent}
for different models in combinatorics and last passage
percolation, in which a few limiting distribution functions were
also computed depending on parameters. But the functions $F_k$,
$k\ge 2$, in Theorem \ref{thm:main} seem to be new in this paper.
The last passage percolation model considered in \cite{BR4,
SpohnPCurrent} has some relevance to our problem; see Section
\ref{sec:last} below.

Theorem \ref{thm:main} and the fact that $F_k$ and $G_k$ are
distribution functions yield the following consequence.

\begin{cor}\label{cor:main}
Under the same assumption of Theorem \ref{thm:main}, the following
holds.
\begin{itemize}
\item[(a)] When for some $0\le k\le r$
\begin{equation}
  \ell_1=\dots= \ell_k=1+\gamma^{-1},
\end{equation}
and $\ell_{k+1}, \dots, \ell_{r}$ are in a compact subset
of $(0, 1+\gamma^{-1})$,
\begin{equation}
  \lambda_1 \to (1+\gamma^{-1})^2 \quad \text{in probability.}
\end{equation}
\item[(b)] When for some $1\le k\le r$,
\begin{equation}
  \ell_1=\cdots=\ell_k > 1+\gamma^{-1}
\end{equation}
and $\ell_{k+1}, \dots, \ell_r$ are in a compact subset of
$(0, \ell_1)$,
\begin{equation}
  \lambda_1 \to
  \ell_1 \biggl(1+\frac{\gamma^{-2}}{\ell_1-1} \biggr)
\quad \text{in probability.}
\end{equation}
\end{itemize}
\end{cor}

\begin{proof}
Suppose we are in the case of (a). For any fixed $\epsilon>0$ and
$x\in\mathbb{R}$,
\begin{equation}
  \limsup_{M\to\infty} \mathbb{P}\bigl( \lambda_1 \le (1-\epsilon) (1+\gamma^{-1})^2 \bigr)
  \le \limsup_{M\to\infty} \mathbb{P}\biggl( \lambda_1\le (1+\gamma^{-1})^2+
  \frac{xM^{1/3}\gamma}{(1+\gamma)^{4/3}} \biggr) =F_k(x).
\end{equation}
By taking $x\to-\infty$, we find that
\begin{equation}
  \lim_{M\to\infty} \mathbb{P}\bigl( \lambda_1 \le (1-\epsilon) (1+\gamma^{-1})^2
  \bigr) =0.
\end{equation}
Similar arguments implies that $\mathbb{P}\bigl( \lambda_1 \ge
(1+\epsilon) (1+\gamma^{-1})^2  \bigr) \to 0.$ The case of (b)
follows from the same argument.
\end{proof}

Together with \eqref{eq:DOS}, Theorem \ref{thm:main}/Corollary
\ref{cor:main} imply that under the Gaussian assumption, when all
but finitely many eigenvalues of $\Sigma$ are $1$, $\lambda_1$ is
separated from the rest of eigenvalues if and only if at least one
eigenvalue of $\Sigma$ is greater than $1+\gamma^{-1}$. Theorem
\ref{thm:main} also claims that when $\lambda_1$ is separated from
the rest, the fluctuation of $\lambda_1$ is of order $M^{1/2}$
rather than $M^{2/3}$. Here the critical value $1+\gamma^{-1}$
comes from a detail of computations and we do not have an
intuitive reason yet. However, see Section \ref{sec:last} below
for a heuristic argument from a last passage percolation model.

Compare the case (b) of Theorem \ref{thm:main}/Corollary
\ref{cor:main} with the following result for samples of
\emph{finite number of variables}.

\begin{prop}\label{prop:finite}
Suppose that there are $M$ samples of $N=k$ variables. Assume that
all the eigenvalues of the covariance matrix are the same;
\begin{equation}
  \ell_1=\cdots=\ell_k.
\end{equation}
Then for fixed $N=k$, as $M\to\infty$,
\begin{equation}
  \lim_{M\to\infty} \mathbb{P} \biggl( \bigl(\lambda_1-\ell_1
  \bigr) \frac1{\ell_1}\sqrt{M} x\biggr) = G_k(x)
\end{equation}
and
\begin{equation}\label{eq:finiteinprob}
  \lambda_1 \to \ell_1 \qquad \text{in probability.}
\end{equation}
\end{prop}

This result shows that the model in the case (b) of Theorem
\ref{thm:main}/Corollary \ref{cor:main} is not entirely dominated
by the distinguished eigenvalues $\ell_1=\cdots=\ell_k$ of the
covariance matrix. Instead the contribution to $\lambda_1$ comes
from both $\ell_1= \cdots=\ell_k$ and infinitely many unit
eigenvalues. The proof of Proposition \ref{prop:finite} is given
in section \ref{sec:finite}.

 Further detailed analysis along the line of this
paper would yield the convergence of the moments of $\lambda_1$
under the scaling of Theorem \ref{thm:main}. This will be
presented somewhere else.

\bigskip

The real question is the real sample covariance. In the null
cases, by comparing \eqref{eq:CLTreal} and \eqref{eq:CLTnull}, we
note that even though the limiting distributions are different,
the scalings are identical. In view of this, we conjecture the
following:

\begin{Conjecture} For real sample covariance, the
Theorem \ref{thm:main} still holds true for different limiting
distributions but \emph{with the same scaling}. In particular, the
critical value of distinguished eigenvalues $\ell_j$ of the
covariance matrix is again expected to be $1+\gamma^{-1}$.
\end{Conjecture}

\subsection{Around the transition point; interpolating distributions}

We also investigate the nature of the transition at
$\ell_j=1+\gamma^{-1}$. The following result shows that if $\ell_j$
themselves scale properly in $M$, there are interpolating limiting
distributions.

We first need more definitions. For $m=1,2,3,\dots$, and for
$w_1,\dots, w_m\in \mathbb{C}$, set
\begin{equation}
  s^{(m)}(u; w_1,\dots, w_m) = \frac1{2\pi}
  \int e^{iua+ i\frac13a^3} \prod_{j=1}^m \frac1{w_j+ia} \, da
\end{equation}
where the contour is from $\infty e^{5i\pi/6}$ to $\infty
e^{i\pi/6}$ such that the points $a=iw_1,\dots, iw_m$ lie
\emph{above} the contour. Also set
\begin{equation}
  t^{(m)}(v; w_1,\dots, w_{m-1}) = \frac1{2\pi} \int e^{ivb+ i\frac13 b^3}
  \prod_{j=1}^{m-1} (w_j-ib) \, db
\end{equation}
where the contour is from $\infty e^{5i\pi/6}$ to $\infty
e^{i\pi/6}$.

\begin{definition}\label{defn:Fkint}
For $k=1,2,\dots$, define for real $x$ and $w_1,\dots, w_k$,
\begin{equation}\label{eq:defFint}
\begin{split}
  & F_k(x;w_1,\dots, w_k) \\
  &\quad = \det \bigl(1-\mathbf{A}_x\bigr)
  \cdot \det\biggl( 1- <\frac1{1-\mathbf{A}_x} s^{(m)}(w_1,\dots, w_m),
  t^{(n)}(w_1,\dots, w_{n-1})> \biggr)_{1\le m,n\le k}.
\end{split}
\end{equation}
\end{definition}

The function $F_1(x;w)$ previously appeared in
\cite{ForresterRinterp} in a disguised form. See (4.18) and (4.12)
of \cite{ForresterRinterp}

The formula \eqref{eq:defFint} may seem to depend on the ordering
of the parameters $w_1,\dots, w_k$.
But as the following result \eqref{eq:main3} shows, it is
independent of the ordering of the parameters. This can also be seen
from a formula of \cite{Baik}.

Like $F_k$, it is not difficult to check that the function
$F_k(x;w_1,\dots, w_k)$ is continuous, non-decreasing and becomes
$1$ as $x\to +\infty$. The proof that $F_k(x;w_1,\dots, w_k) \to
0$ as $x\to -\infty$ is in \cite{Baik}. Therefore,
$F_k(x;w_1,\dots, w_k)$ is a distribution function. It is direct
to check that $F_k(x;w_1,\dots, w_k)$ interpolates $F_0(x), \dots,
F_k(x)$. For example, $\lim_{w_2\to +\infty} F_2(x, 0, w_2) =
F_1(x)$, $\lim_{w_1\to +\infty} \lim_{w_2\to +\infty}
F_2(x;w_1,w_2) = F_0(x)$, etc.

\begin{theo}\label{thm:main3}
Suppose that for a fixed $r$, $\ell_{r+1}=\ell_{r+2}=\dots
=\ell_N=1$. Set for some $1\le k\le r$,
\begin{equation}
  \ell_j = 1+\gamma^{-1} - \frac{(1+\gamma)^{2/3}w_j}{\gamma M^{1/3}},
  \qquad j=1,2,\dots, k.
\end{equation}
When $w_j$, $1\le j\le k$, is in a compact subset of $\mathbb{R}$,
and $\ell_{j}$, $k+1\le j\le r$, is in a compact subset of
$(0,1+\gamma^{-1})$, as $M, N\to\infty$ such that $M/N=\gamma^2$
is in a compact subset of $[1,\infty)$,
\begin{equation}\label{eq:main3}
\mathbb{P} \biggl( \bigr( \lambda_1 - (1+\gamma^{-1})^2 \bigr)
\cdot
  \frac{\gamma}{(1+\gamma)^{4/3}}M^{2/3} \le x \biggr)
  \to F_k(x;w_1,\dots, w_k)
\end{equation}
for any $x$ in a compact subset of $\mathbb{R}$.
\end{theo}

The Painlev\'e II-type expression for $F_k(x;w_1,\dots, w_k)$
will be presented in \cite{Baik}.



\bigskip

This paper is organized as follows. The basic algebraic formula of
the distribution of $\lambda_1$ in terms of a Fredholm determinant
is given in Section \ref{sec:basicformulas}, where an outline of
the asymptotic analysis of the Fredholm determinant is also
presented. The proofs of Theorem \ref{thm:main} (a) and Theorem
\ref{thm:main3} are given in Section \ref{sec:Airylimit}. The
proof of Theorem \ref{thm:main} (b) is in Section
\ref{sec:Gausslimit} and the proof of Proposition
\ref{prop:finite} is presented in Section \ref{sec:finite}. In
Section \ref{sec:last}, we indicate a connection between the
sample covariance matrices, and a last passage percolation model
and also a queueing theory.

\begin{nrem}
Throughout the paper, we set
\begin{equation}
  \pi_j=\ell_j^{-1}.
\end{equation}
This is only because the formulas below involving $\ell_j^{-1}$
become simpler with $\pi_j$.
\end{nrem}

\bigskip

\noindent {\bf Acknowledgments.} Special thanks is due to Iain
Johnstone for kindly explaining the importance of computing the
largest eigenvalue distribution for non-null covariance case and
also for his constant interest and encouragement. We would like to
thank Kurt Johansson for sharing with us his  proof of Proposition
\ref{prop:integralformula} below, Eric Rains for useful
discussions and also Debashis Paul for finding a typographical
error in the main theorem in an earlier draft. The work of the
first author was supported in part by NSF Grant \# DMS-0350729.

\section{Basic formulas}\label{sec:basicformulas}

\begin{nrem}
The notation $V(x)$ denotes the
Vandermonde determinant
\begin{equation}
  V(x)=\prod_{i<j} (x_i-x_j)
\end{equation}
of a (finite) sequence $x=(x_1,x_2,\dots)$.
\end{nrem}

\subsection{Eigenvalue density; algebraic formula}

For complex Gaussian samples, the density of the sample covariance
matrix is already known to Wishart around 1928 (see e.g.
\cite{Muirhead}):
\begin{equation}
  p(S) = \frac{1}{C} e^{-M \cdot tr(\Sigma^{-1} S)} (\det
  S)^{M-N}
\end{equation}
for some normaliztion constant $C>0$. As $\Sigma$ and $S$ are
Hermitian, we can set $\Sigma=UDU^{-1}$ and $S=HLH^{-1}$ where $U$
and $H$ are unitary matrices, $D= diag(\ell_1, \cdots, \ell_N)
=diag(\pi_1^{-1},\cdots, \pi_N^{-1})$ and $L =
diag(\lambda_1,\cdots, \lambda_N)$. By taking the change of
variables $S \mapsto (L,H)$ using the Jacobian formula $dS=
cV(L)^2dLdH$ for some constant $c>0$, and then integrating over
$H$, the density of the eigenvalues is (see, e.g.
\cite{JamesAnn35})
\begin{equation}\label{eq:eigendenoverU}
  p(\lambda) = \frac1{C} V(\lambda)^2 \prod_{j=1}^N
  \lambda_j^{M-N} \cdot
  \int_{Q\in U(N)} e^{-M\cdot tr(D^{-1}QLQ^{-1})} dQ.
\end{equation}
for some (new) constant $C>0$ where $U(N)$ is the set of $N\times
N$ unitary matrices and $\lambda=(\lambda_1,\dots, \lambda_N)$.
The last integral is known as Harish-Chandra-Itzykson-Zuber
integral (see e.g. \cite{Mehta}) and we find
\begin{equation}\label{eq:eigenden}
  p(\lambda) = \frac1{C}
  \frac{\det(e^{-M\pi_j\lambda_k})_{1\le j,k\le N}}{V(\pi)}
  V(\lambda)
  \prod_{j=1}^N  \lambda_j^{M-N}.
\end{equation}
Here when some of $\pi_j$'s coincide, we interpret the formula
using the l'Hopital's rule. We note that for a real sample
covariance matrix, it is not known if the corresponding integral
over the orthogonal group $O(N)$ is computable as above. Instead
one usually define hypergeometric functions of matrix argument and
study their algebraic properties (see e.g. \cite{Muirhead}).
Consequently, the techniques below that we will use for the
density of the form \eqref{eq:eigenden} is not applicable to real
sample matrices.

For the density \eqref{eq:eigenden}, the distribution function of
the largest eigenvalue $\lambda_1$  can be expressed in terms of a
Fredholm determinant, which will be the starting point of our
asymptotic analysis. The following result can be obtained by
suitably re-interpreting and taking a limit of a result of
\cite{Okounkov}. A different proof is given in \cite{Sandrine}.
For the convenience of reader we include yet another proof by
Johansson \cite{Joh:unp} which uses an idea from random matrix
theory (see e.g. \cite{TW2}).

\begin{figure}[ht]
\centerline{\epsfxsize=8cm\epsfbox{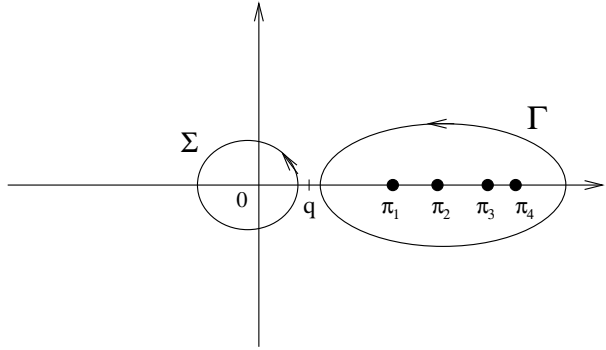}}
\caption{Contours $\Gamma$ and $\Sigma$}\label{fig:Kcontour}
\end{figure}

\begin{prop}\label{prop:integralformula}
For any fixed $\qq$ satisfying $0<\qq<\min\{\pi_j\}_{j=1}^N$,
let $\mathbf{K}_{M,N}|_{(\xi,\infty)}$ be the operator acting on
$L^2((\xi,\infty))$ with kernel
\begin{equation}\label{eq:Kkernel}
  \mathbf{K}_{M,N}(\eta, \zeta) =
  \frac{M}{(2\pi i)^2} \int_\Gamma dz \int_\Sigma dw
  \,
  e^{- \eta M(z-\qq) + \zeta M(w-\qq)}
  \frac1{w-z} \bigl(\frac{z}{w}\bigr)^M
  \prod_{k=1}^N
  \frac{\pi_k-w}{\pi_k-z}
\end{equation}
where $\Sigma$ is a simple closed contour enclosing $0$ and lying
in $\{w:Re(w)<\qq\}$, and $\Gamma$ is a simple closed contour
enclosing $\pi_1,\dots, \pi_N$ and lying $\{z:Re(z)>\qq\}$, both
oriented counter-clockwise (see Figure \ref{fig:Kcontour}). Then
for any $\xi\in\mathbb{R}$,
\begin{equation}\label{eq:Fred}
  \mathbb{P} (\lambda_1 \le \xi)
  = \det( 1-\mathbf{K}_{M,N}|_{(\xi,\infty)}).
\end{equation}
\end{prop}

\begin{remark}
Note that the left-hand-side of \eqref{eq:Fred} does not depend on
the parameter $q$. The Fredholm determinant on the right-hand-side
of \eqref{eq:Fred} is also independent of the choice of $\qq$ as
long as $0<\qq<\min\{\pi_j\}_{j=1}^N$. If we use the notation
$\mathbf{K}_q$ to denote $\mathbf{K}_{M,N}|_{(\xi,\infty)}$ for
the parameter $q$, then
$\mathbf{K}_{q'}=\mathbf{E}\mathbf{K}_q\mathbf{E}^{-1}$ where
$\mathbf{E}$ is the multiplication by $e^{(\qq'-\qq)\lambda}$;
$(\mathbf{E}f)(\lambda)= e^{(\qq'-\qq)\lambda} f(\lambda)$. But
determinants are invariant under conjugations as long as both
$\mathbf{K}_{q'}$ and $\mathbf{E}\mathbf{K}_q\mathbf{E}^{-1}$ are
in the trace class, which is the case when $0<\qq, \qq'
<\min\{\pi_j\}_{j=1}^N$. The parameter $\qq$ ensures that the
kernel $\mathbf{K}_{M,N}(\eta,\zeta)$ is finite when $\eta\to
+\infty$ or $\zeta \to +\infty$ and the operator
$\mathbf{K}_{M,N}|_{(\xi,\infty)}$ is trace class. It also helps
the proof of the convergence of the operator in the next section.
\end{remark}

\begin{proof}
For a moment we assume that all $\pi_j$'s are distinct. Note that
the density \eqref{eq:eigenden} is symmetric in $\lambda_j$'s.
Hence using $V(\lambda)=\det(\lambda_k^{j-1})$, we find that
\begin{equation}
  \mathbb{P}(\lambda_1\le \xi)
  = \frac1{C'} \int_0^\infty \cdots \int_0^\infty
  \det(\lambda_k^{j-1})
  \det(e^{-M\pi_j \lambda_k})  \prod_{k=1}^N
  (1-\chi_{(\xi,\infty)}(\lambda_k))\lambda_k^{M-N} d\lambda_k
\end{equation}
with some constant $C'>0$, where $\chi_{(\xi,\infty)}$ denotes the
characteristic function (indicator function). Using the
fundamental identity which dates back to \cite{Andreief},
\begin{equation}
  \int\cdots \int \det(f_j(x_k)) \det(g_j(x_k)) \prod_k d\mu(x_k)
  = \det \biggl( \int f_j(x)g_k(x) d\mu(x) \biggr),
\end{equation}
we find
\begin{equation}
  \mathbb{P} (\lambda_1\le \xi)
  = \frac1{C'} \det \biggl(
  \int_0^\infty (1-\chi_{(\xi,\infty)}(\lambda)) \lambda^{j-1+M-N} e^{-M\pi_k \lambda}
  d\lambda
  \biggr)_{1\le j,k\le N} .
\end{equation}
Now set $\nu=M-N$,
$\phi_j(\lambda)=\lambda^{j-1+\nu}e^{-M\qq\lambda}$ and
$\Phi_k(\lambda)=e^{-M(\pi_k-\qq) \lambda}$ for any $\qq$
such that $0<\qq<\min\{\pi_j\}_{j=1}^N$. Also let
\begin{equation}
  A=(A_{jk})_{1\le j,k\le N}, \qquad A_{jk} = \int_0^{\infty} \phi_j(\lambda)
  \Phi_k(\lambda)
  d\lambda
  = \frac{\Gamma(j+\nu)}{(M\pi_k)^{j+\nu}}.
\end{equation}
A direct computation shows that
\begin{eqnarray}\label{eq:detA}
 \det A
 = \prod_{j=1}^N \frac{\Gamma(j+\nu)}{(M\pi_j)^{\nu+1}}
 \cdot \det((M\pi_j)^{-(j-1)})
 = \prod_{j=1}^N \frac{\Gamma(j+\nu)}{(M\pi_j)^{\nu+1}}
 \cdot \prod_{1\le j<k\le N} ((M\pi_j)^{-1}-(M\pi_k)^{-1}).
\end{eqnarray}
Thus $A$ is invertible. Also define the operators
$B:L^2((0,\infty))\to \ell^2(\{1,\dots, N\})$, $C:
\ell^2(\{1,\dots, N\}) \to L^2((0,\infty))$ by
\begin{equation}
  B(j,\lambda) = \phi_j(\lambda), \qquad C(\lambda, k) =
  \Phi_k(\lambda)
\end{equation}
and let $P_\xi$ be the projection from $(0,\infty)$ to
$(\xi,\infty)$. Then as
\begin{equation}
  \int_0^\infty \chi_{(\xi,\infty)}(\lambda) \lambda^{j-1+M-N} e^{-M\pi_k \lambda}
  d\lambda
  = (BP_\xi C)(j,k),
\end{equation}
we find that
\begin{equation}
  \mathbb{P}(\lambda_1\le \xi)
  = \frac1{C'} \det(A-BP_\xi C).
\end{equation}
So,
\begin{equation}\label{eq:detalgt1}
  \mathbb{P}(\lambda_1\le \xi)
  = \frac{\det(A)}{C'} \det(1-A^{-1}B P_\xi C)
  = C'' \det(1-P_\xi CA^{-1}B)
  = C'' \det(1-P_\xi CA^{-1}B P_\xi)
\end{equation}
for some constant $C''$ which does not depend on $\xi$. But by
letting $\xi\to +\infty$ in both sides of \eqref{eq:detalgt1}, we
easily find that $C''=1$. The kernel of the operator in the
determinant is
\begin{equation}
  (CA^{-1}B)(\eta,\zeta) =  \sum_{j=1}^N C(\eta, j)
  (A^{-1}B)(j,\zeta), \qquad \eta, \zeta>\xi,
\end{equation}
and from the Cramer's rule,
\begin{equation}
  (A^{-1}B)(j,\zeta) = \frac{\det A^{(j)}(\zeta)}{\det A}
\end{equation}
where $A^{(j)}(\zeta)$ is the matrix given by $A$ with $j$th
column replaced by the vector $(\phi_1(\zeta),\cdots,
\phi_N(\zeta))^T$. To compute $A^{(j)}$, note (the Hankel's
formula for Gamma function) that for a positive integer $a$
\begin{equation}
  \frac{1}{2\pi i}\int_{\Sigma} \frac{e^w}{w^a} dw
  = \frac1{(a-1)!} = \frac{1}{\Gamma(a)}
\end{equation}
where $\Sigma$ is any simple closed contour enclosing the origin
$0$ with counter-clockwise orientation.
By replacing $w\to \zeta Mw$ and setting $a=j+\nu$, this implies
that
\begin{equation}
  \zeta^{j-1+\nu}=  \frac{\Gamma(j+\nu)}{2\pi i
}\int_{\Sigma}e^{\zeta Mw}\frac{M}{(Mw)^{j+\nu}}dw.
\end{equation}
Substituting this formula for $\phi_j(\zeta)$ in the $j$th column
of $A^{(j)}$, and pulling out the integrals over $w$,
\begin{equation}
  \det A^{(j)}(\zeta) = \frac1{2\pi i} \int_\Sigma e^{\zeta M(w-\qq)}
  \det( A'(w)) Mdw
\end{equation}
where the entries of $A'(w)$ are $A'_{ab}(w)=
\Gamma(a+\nu)/p_{b}^{a+\nu}$ where $p_b=M\pi_b$ when $b\neq j$ and
$p_b=Mw$ when $b=j$. Hence
\begin{equation}
  \det A^{(j)}(\zeta)
  = \prod_{k\neq j} \frac1{(M\pi_k)^{1+\nu}} \cdot
  \prod_{k=1}^N \Gamma(k+\nu) \cdot
  \frac1{2\pi i} \int_\Sigma e^{\zeta M(w-\qq)} \prod_{1\le a<b \le N}(p_a^{-1}-p_b^{-1})
   \frac{Mdw}{(Mw)^{1+\nu}},
\end{equation}
and so using \eqref{eq:detA},
\begin{equation}
  (A^{-1}B)(j,\zeta) =
  \frac{M\pi_j^{N+\nu}}{2\pi i} \int_\Sigma e^{\zeta M(w-\qq)}
  \prod_{k\neq j} \frac{w-\pi_k}{\pi_j-\pi_k}
  \frac{dw}{w^{N+\nu}} .
\end{equation}
But for any simple closed contour $\Gamma$ that encloses
$\pi_1,\dots, \pi_N$ but excludes $w$, and is oriented
counter-clockwise,
\begin{equation}
  \frac1{2\pi i} \int_\Gamma  z^M e^{-\eta Mz} \frac1{w-z}
  \prod_{k=1}^N \frac{w-\pi_k}{z-\pi_k} \, dz
  = \sum_{j=1}^N \pi_j^M e^{-M\pi_j \eta} \prod_{k\neq j}
  \frac{w-\pi_k}{\pi_j-\pi_k}.
\end{equation}
Therefore, we find (note $N+\nu=M$)
\begin{equation}
  (CA^{-1}B)(\eta,\zeta)
  = \frac{M}{(2\pi i)^2} \int_\Gamma dz \int_\Sigma dw
  \,
  e^{- \eta M(z-\qq) + \zeta M(w-\qq)} \frac1{w-z}
  \prod_{k=1}^N
  \frac{w-\pi_k}{z-\pi_k} \cdot \bigl(\frac{z}{w}\bigr)^M
\end{equation}
which completes the proof when all $\pi_j$'s are distinct. When
some of $\pi_j$'s are identical, the result follows by taking
proper limits and using the l'Hospital's theorem.
\end{proof}

Note that for $z\in\Gamma$ and $w\in\Sigma$, $Re(w-z)<0$. Hence
using
\begin{equation}\label{eq:63}
  \frac{1}{w-z} = -M \int_0^\infty e^{yM(w-\qq-(z-\qq))}dy
\end{equation}
for $1/(w-z)$ in \eqref{eq:Kkernel}, the kernel
$\mathbf{K}_{M,N}(\eta,\zeta)$ is equal to
\begin{equation}
  \mathbf{K}_{M,N}(\eta,\zeta) =\int_0^\infty \mathbf{H}(\eta+y)
  \mathbf{J}(\zeta+y) dy
\end{equation}
where
\begin{equation}\label{eq:64}
  \mathbf{H}(\eta+y) = \frac{M}{2\pi} \int_\Gamma e^{-(\eta+y)M(z-\qq)}
  z^M
  \prod_{k=1}^N \frac1{\pi_k-z} dz
\end{equation}
and
\begin{equation}\label{eq:65}
  \mathbf{J}(\zeta+y) = \frac{M}{2\pi} \int_\Sigma
  e^{(\zeta+y)M(w-\qq)}w^{-M}
  \prod_{k=1}^N (\pi_k-w) dw.
\end{equation}

\subsection{Asymptotic analysis: basic ideas}\label{sec:basicas}

From now on, as mentioned in the Introduction, we assume that
\begin{equation}
  \pi_{r+1}=\dots = \pi_N=1.
\end{equation}
In this case, \eqref{eq:63}, \eqref{eq:64} and \eqref{eq:65} become
\begin{equation}\label{eq:KinHJ}
  \mathbf{K}_{M,N}(\eta,\zeta) =\int_0^\infty \mathbf{H}(\eta+y)
  \mathbf{J}(\zeta+y) dy
\end{equation}
where
\begin{equation}
  \mathbf{H}(\eta) = \frac{M}{2\pi} \int_\Gamma e^{-M\eta (z-\qq)}
  \frac{z^M}{(1-z)^{N-r}}
  \prod_{k=1}^r \frac1{\pi_k-z} dz
\end{equation}
and
\begin{equation}
  \mathbf{J}(\zeta) = \frac{M}{2\pi} \int_\Sigma
  e^{M\zeta (z-\qq)}\frac{(1-z)^{N-r}}{z^M}
  \prod_{k=1}^r (\pi_k-z) dz.
\end{equation}

Set
\begin{equation}\label{eq:MN}
  \frac{M}{N}=\gamma^2\ge 1.
\end{equation}
For various choices of $\pi_j$, $1\le j\le r$, we will consider
the limit of $\mathbb{P}(\lambda_1\le \xi)$ when $\xi$ is scaled
as of the form (see Theorem \ref{thm:main})
\begin{equation}\label{eq:scalegeneral}
  \xi= \mu+\frac{\nu x}{M^\alpha}
\end{equation}
for some constants $\mu=\mu(\gamma), \nu=\nu(\gamma)$ and for some
$\alpha$, while $x$ is a fixed real number. By translation and
scaling, the equation \eqref{eq:Fred} becomes
\begin{equation}
  \mathbb{P}\biggl(\lambda_1 \le \mu+\frac{\nu x}{M^\alpha}\biggr)
  = \det( 1-\mathbf{K}_{M,N}|_{(\mu+\frac{\nu x}{M^\alpha}, \infty)})
  = \det(1 - \mathcal{K}_{M,N})
\end{equation}
where $\mathcal{K}_{M,N}$ is the operator acting on
$L^2((0,\infty))$ with kernel
\begin{equation}
  \mathcal{K}_{M,N}(u,v) = \frac{\nu}{M^{\alpha}}
  \mathbf{K}_{M,N} \biggl(\mu+ \frac{\nu (x+u)}{M^\alpha},
  \mu+ \frac{\nu (x+u)}{M^\alpha}\biggr).
\end{equation}
Using \eqref{eq:KinHJ}, this kernel is equal to
\begin{equation}\label{eq:KcalinHJ}
  \mathcal{K}_{M,N}(u,v) =\int_0^\infty \mathcal{H}(x+u+y)
  \mathcal{J}(x+v+y) dy
\end{equation}
where
\begin{equation}\label{eq:Hcal}
  \mathcal{H}(u) = \frac{\nu M^{1-\alpha}}{2\pi}
  \int_\Gamma e^{-\nu M^{1-\alpha} u(z-\qq)} e^{-M\mu (z-\qq)}
  \frac{z^M}{(1-z)^{N-r}}
  \prod_{\ell=1}^r \frac1{\pi_\ell-z} dz
\end{equation}
and
\begin{equation}\label{eq:Jcal}
  \mathcal{J}(v) = \frac{\nu M^{1-\alpha}}{2\pi} \int_{\Sigma}
  e^{\nu M^{1-\alpha} v(w-\qq)} e^{M\mu (w-\qq)}\frac{(1-w)^{N-r}}{w^M}
  \prod_{\ell=1}^r (\pi_\ell-w) dw.
\end{equation}

We need to find limits of $\mathcal{K}_{M,N}(u,v)$ for various
choices of $\pi_j$'s as $M,N\to\infty$. A sufficient condition for
the convergence of a Fredholm determinant is the convergence in
trace norm of the operator. As $\mathcal{K}_{M,N}$ is a product of
two operators, it is enough to prove the convergences of
$\mathcal{H}$ and $\mathcal{J}$ in Hilbert-Schmidt norm. Hence in
sections \ref{sec:Airylimit} and \ref{sec:Gausslimit} below, we
will prove that for proper choices of $\mu, \nu$ and $\alpha$,
there are limiting operators $\mathcal{H}_\infty$ and
$\mathcal{J}_\infty$ acting on $L^2((0,\infty))$ such that for any
real $x$ in a compact set,
\begin{equation}\label{eq:HS1}
  \int_0^\infty \int_0^\infty
  \bigl| Z_M \mathcal{H}_{M,N}(x+u+y)-\mathcal{H}_\infty(x+u+y) \bigr|^2 dudy
  \to 0
\end{equation}
and
\begin{equation}\label{eq:HS2}
  \int_0^\infty \int_0^\infty
  \bigl| \frac1{Z_M} \mathcal{J}_{M,N}(x+u+y)-\mathcal{J}_\infty(x+u+y) \bigr|^2 dudy
  \to 0
\end{equation}
for some non-zero constant $Z_M$ as $M, N \to\infty$ satisfying
\eqref{eq:MN} for $\gamma$ in a compact set. We will use
steepest-descent analysis.

\section{Proofs of Theorem \ref{thm:main} (a) and Theorem \ref{thm:main3}}
\label{sec:Airylimit}

We first consider the proof of Theorem \ref{thm:main} (a). The
proof of Theorem \ref{thm:main3} will be very similar (see the
subsection \ref{sec:main3} below).
We assume that for some $0\le k\le r$,
\begin{equation}
  \pi_1^{-1}=\dots= \pi_k^{-1}=1+\gamma^{-1}
\end{equation}
 and
$\pi_{k+1}^{-1}\ge\dots\ge \pi_{r}^{-1}$ are in a compact subset
of $(0, 1+\gamma^{-1})$.

For the scaling \eqref{eq:scalegeneral}, we take
\begin{equation}
\alpha=2/3
\end{equation}
and
\begin{equation}\label{eq:mu}
  \mu=\mu(\gamma) := \biggl(\frac{1+\gamma}{\gamma}\biggr)^2, \qquad
  \nu=\nu(\gamma) := \frac{(1+\gamma)^{4/3}}{\gamma}
\end{equation}
so that
\begin{equation}
  \xi= \mu+\frac{\nu x}{M^{2/3}}.
\end{equation}
The reason for such choices will be made clear during the
following asymptotic analysis. There is still an arbitrary
parameter $\qq$. It will be chosen in \eqref{eq:c1} below.

The functions \eqref{eq:Hcal} and \eqref{eq:Jcal} are now
\begin{equation}\label{eq:H1}
  \mathcal{H}(u)
  = \frac{\nu M^{1/3}}{2\pi}
  \int_\Gamma e^{-\nu M^{1/3} u(z-\qq)} e^{Mf(z)}
  \frac1{\bigl(p_c-z\bigr)^kg(z)} dz
\end{equation}
and
\begin{equation}\label{eq:J1}
  \mathcal{J}(v) = \frac{\nu M^{1/3}}{2\pi} \int_\Sigma
  e^{\nu M^{1/3} v(z-\qq)} e^{-Mf(z)}
  \bigl(p_c-z\bigr)^kg(z) dz
\end{equation}
with
\begin{equation}
  p_c := \frac{\gamma}{\gamma+1},
\end{equation}
and
\begin{equation}\label{eq:f}
   f(z) := -\mu (z-\qq) +\log(z) - \frac1{\gamma^2}\log(1-z),
\end{equation}
where we take the principal branch of $\log$ (i.e. $\log(z)= \ln |z|+i
arg(z)$, $-\pi<arg(z)<\pi$), and
\begin{equation}\label{eq:g}
  g(z) := \frac{1}{(1-z)^{r}}
  \prod_{\ell=k+1}^r (\pi_\ell-z).
\end{equation}

Now we find the critical point of $f$. As
\begin{equation}
  f'(z) = -\mu + \frac{1}{z} -\frac1{\gamma^2(z-1)},
\end{equation}
$f'(z)=0$ is a quadratic equation. But with the choice
\eqref{eq:mu} of $\mu$, there is a double root at
$z=p_c=\frac{\gamma}{\gamma+1}$. Note that for $\gamma$ in a
compact subset of $[1,\infty)$, $p_c$ is strictly less than $1$.
Being a double root,
\begin{equation}
  f'(p_c)=f''(p_c)=0,
\end{equation}
where
\begin{equation}
  f''(z) = - \frac{1}{z^2} + \frac1{\gamma^2(z-1)^2}.
\end{equation}
It is also direct to compute
\begin{equation}\label{eq:fthreep}
  f^{(3)}(z) = \frac{2}{z^3} - \frac{2}{\gamma^2(z-1)^3},
  \qquad f^{(3)}(p_c) = \frac{2(\gamma+1)^4}{\gamma^3} = 2\nu^3
\end{equation}
and
\begin{equation}
  f(p_c) = -\mu(p_c-\qq) +
  \log\biggl(\frac{\gamma}{\gamma+1}\biggr) -
  \frac1{\gamma^2} \log\biggl( \frac1{\gamma+1} \biggr).
\end{equation}

\begin{figure}[ht]
\centerline{\epsfxsize=17cm\epsfbox{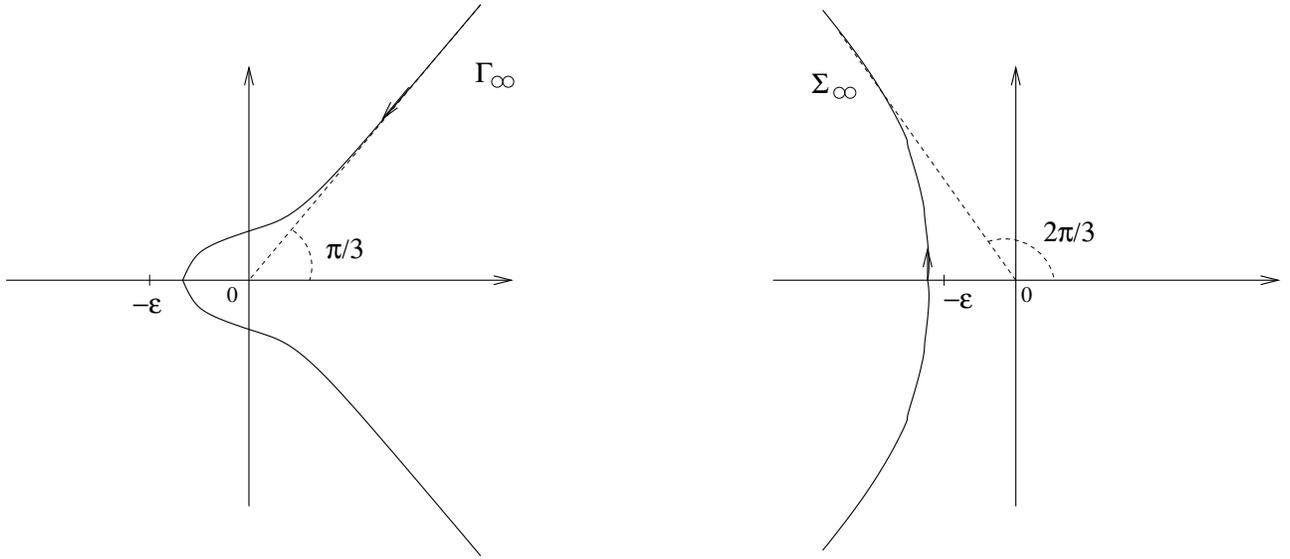}}
\caption{Contours $\Gamma_\infty$ and $\Sigma_\infty$
of $\mathcal{H}(u)$ and $\mathcal{J}(v)$}\label{fig:HJcon}
\end{figure}

As $f^{(3)}(p_c)>0$, the steepest-descent curve of $f(z)$ comes to
the point $p_c$ with angle $\pm \pi/3$ to the real axis. Once the
contour $\Gamma$ is chosen to be the steepest-descent curve near
$p_c$ and is extended properly, it is expected that the main
contribution to the integral of $\mathcal{H}$ comes from a contour
near $p_c$. There is, however, a difficulty since at the critical
point $z=p_c$, the integral \eqref{eq:H1} blows up due to the term
$(p_c-z)^k$ in the denominator when $k\ge 1$. Nevertheless this
can be overcome if we choose $\Gamma$ to be close to $p_c$, but
not exactly pass through $z_0$. Also as the contour $\Gamma$
should contain all $\pi_j$'s, some of which may be equal to $p_c$,
we will choose $\Gamma$ to intersect the real axis \emph{to the
left} of $p_c$. By formally approximating the function $f$ by a
third-degree polynomial and the function $g(z)$ by $g(p_c)$, we
expect
\begin{equation}
\begin{split}
  \mathcal{H}(u)
  & \sim
  \frac{\nu M^{1/3}}{2\pi g(p_c)}
  \int_\Gamma  e^{-\nu M^{1/3} u(z-\qq)} e^{M\bigl(f(z_c)+\frac{f^{(3)}(p_c)}{3!}(z-p_c)^3\bigr)}
  \frac1{\bigl(p_c-z\bigr)^k} dz
\end{split}
\end{equation}
for some contour $\Gamma_\infty$. Now taking the intersection
point of $\Gamma$ with the real axis to be on the left of $p_c$ of
distance of order $M^{-1/3}$, and then changing of the variables
by $\nu M^{1/3} (z-p_c)=a$, we expect
\begin{equation}
\begin{split}
  \mathcal{H}(u)
  \sim
  \frac{(-\nu M^{1/3})^k  e^{Mf(p_c)}}{2\pi g(p_c)}
  e^{-\nu M^{1/3} u(p_c-\qq)} \int_{\Gamma_\infty}  e^{-ua+\frac13a^3}
  \frac1{a^k} da
\end{split}
\end{equation}
Similarly, we expect that
\begin{equation}
\begin{split}
  \mathcal{J}(v) \sim
  \frac{g(p_c) e^{-Mf(p_c)}}{2\pi (-\nu M^{1/3})^k}
  e^{\nu M^{1/3} v(p_c-\qq)} \int_{\Sigma_\infty}  e^{va-\frac13a^3} a^k da
\end{split}
\end{equation}
for some contour $\Sigma_\infty$. When multiplying
$\mathcal{H}(u)$ and $\mathcal{J}(v)$, the constant prefactors
\begin{equation}\label{eq:C1}
  Z_M := \frac{g(p_c)}{(-\nu M^{1/3})^k  e^{Mf(p_c)}}
\end{equation}
and $1/Z_M$ cancel each other out. However, note that there are
still functions $e^{-\nu M^{1/3} u(p_c-\qq)}$ and $e^{\nu M^{1/3}
v(p_c-\qq)}$, one of which may become large when $M\to\infty$ (cf.
\eqref{eq:HS1}, \eqref{eq:HS2}). This trouble can be avoided if we
can simply take $\qq=p_c$, but since $\Gamma$ should be on the
left of $\qq$, this simple choice is excluded. Nevertheless, we
can still take $\qq$ to be $p_c$ minus some positive constant of
order $M^{-1/3}$. Fix
\begin{equation}
 \epsilon>0
\end{equation}
and set
\begin{equation}\label{eq:c1}
  \qq :=p_c - \frac{\epsilon}{\nu M^{1/3}}.
\end{equation}
We then take $\Gamma$ to intersect the real axis at $p_c - c/(\nu
M^{1/3})$ for some $0<c<\epsilon$. With this choice of $\qq$ and
$\Gamma, \Sigma$, we expect that
\begin{equation}
  Z_M \mathcal{H}(u) \sim
  \mathcal{H}_\infty(u)
\end{equation}
where
\begin{equation}\label{eq:Hinf1}
  \mathcal{H}_\infty(u) := \frac{e^{-\epsilon u}}{2\pi}
\int_{\Gamma_{\infty}}  e^{-ua+\frac13a^3}
  \frac1{a^k} da,
\end{equation}
and
\begin{equation}
   \frac1{Z_M} \mathcal{J}(v) \sim \mathcal{J}_\infty(v)
\end{equation}
where
\begin{equation}\label{eq:Jinf1}
  \mathcal{J}_\infty(v) := \frac{e^{\epsilon v}}{2\pi}  \int_{\Sigma_\infty}
e^{va-\frac13a^3}
  a^k da.
\end{equation}
Here the contour $\Gamma_\infty$ is, as in the left picture of
Figure \ref{fig:HJcon}, from $\infty e^{i\pi/3}$ to $\infty
e^{-i\pi/3}$, passes the real axis on the left of the origin and
lies in the region $Re(a+\epsilon)>0$, is symmetric about the real
axis and is oriented from top to bottom. The contour
$\Sigma_\infty$ is, as in the right picture of Figure
\ref{fig:HJcon}, from $\infty e^{-i2\pi/3}$ to $\infty
e^{i2\pi/3}$, lies on the region $Re(a+\epsilon)<0$, is symmetric
about the real axis and is oriented from bottom to top.

This argument should be justified and the following is a rigorous
estimate.

\begin{prop}\label{prop:HJmaina}
Fix $\epsilon>0$ and set $\qq$ by \eqref{eq:c1}. Define
$\mathcal{H}_\infty(u)$ and $\mathcal{J}_\infty(v)$ by
\eqref{eq:Hinf1} and \eqref{eq:Jinf1}, respectively. Then the
followings hold for $\mathcal{H}(u)$ and $\mathcal{J}(v)$ in
\eqref{eq:H1} and \eqref{eq:J1} for $M/N=\gamma^2$ with $\gamma$
in a compact subset of $[1,\infty)$.
\begin{itemize}
\item[(i)] For any fixed $U\in\mathbb{R}$,
there are constants $C, c>0, M_0>0$ such that
\begin{equation}\label{eq:estprop1}
  \bigl| Z_M \mathcal{H}(u) - \mathcal{H}_\infty(u) \bigr|
  \le \frac{Ce^{-cu}}{M^{1/3}}
\end{equation}
for $u\ge U$ when $M\ge M_0$.
\item[(ii)] For any fixed $V\in \mathbb{R}$,
there are constants $C,c>0, M_0>0$ such that
\begin{equation}\label{eq:estprop3}
  \bigl| \frac1{Z_M}\mathcal{J}(v)
- \mathcal{J}_\infty(v) \bigr|
  \le \frac{Ce^{-cv}}{M^{1/3}}
\end{equation}
for $v\ge V$ when $M\ge M_0$.
\end{itemize}
\end{prop}

The proof of this result is given in the following two subsections.

\bigskip

\subsection{Proof of Proposition \ref{prop:HJmaina} (i)}

\begin{figure}[ht]
\centerline{\epsfxsize=11cm\epsfbox{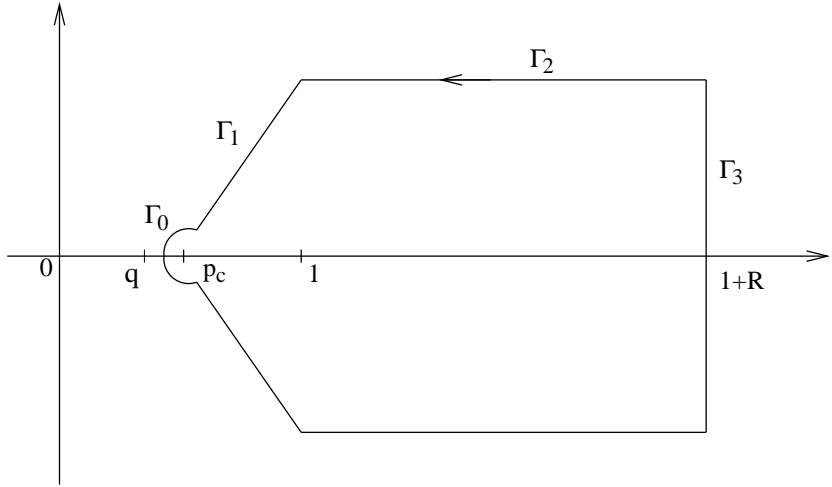}} \caption{Contour
$\Gamma$}\label{fig:Steep1}
\end{figure}

The steepest-descent curve of $f$ will depend on $\gamma$ and $M$.
Instead of controlling uniformity of the curve in $\gamma$ and
$M$, we will rather explicitly choose $\Gamma$ which will be a
steep-descent (though not the steepest-descent) curve of $f(z)$.
Fix $R>0$ such that $1+R>\max\{1,\pi_{r+1},\dots, \pi_k\}$. Define
\begin{eqnarray}
  \Gamma_0 &:=& \{ p_c+ \frac{\epsilon}{2\nu M^{1/3}}e^{i\theta} :
\pi/3\le \theta\le
  \pi\} \\
  \Gamma_1 &:=& \{ p_c+ te^{i\pi/3} :
\frac{\epsilon}{2\nu M^{1/3}} \le t \le 2(1-p_c) \} \\
  \Gamma_2 &:=& \{ p_c+2(1-p_c)e^{i\pi/3} + x : 0\le x\le R\} \\
  \Gamma_3 &:=& \{ 1+R+iy : 0\le y\le \sqrt{3}(1-p_c) \}.
\end{eqnarray}
Set
\begin{equation}
  \Gamma = \bigl( \cup_{k=0}^3 \Gamma_k \bigr) \cup \overline{\bigl( \cup_{k=0}^3 \Gamma_k \bigr)}
\end{equation}
and choose the orientation of $\Gamma$ counter-clockwise. See
Figure \ref{fig:Steep1}.

Note that all the singular points of the integrand of
$\mathcal{H}$ are inside of $\Gamma$ and hence the deformation to
this new $\Gamma$ is allowed. Direct calculations show the
following properties of $\Gamma$. Recall \eqref{eq:f},
\begin{equation}
   f(z) := -\mu (z-\qq) +\log(z) - \frac1{\gamma^2}\log(1-z).
\end{equation}

\begin{lemma}\label{lem:Refdecreasing}
For $\gamma\ge 1$, $Re(f(z))$ is decreasing for
$z\in\Gamma_1\cup\Gamma_2$ as $Re(z)$ increases. Also for $\gamma$
in a compact subset of $[1,\infty)$, we can take $R>0$ large enough
such that
\begin{equation}
  \max_{z\in\Gamma_3} Re(f(z)) \le Re \bigl(f(p_*)\bigr).
\end{equation}
where $p_*:=p_c+2(1-p_c)e^{i\pi/3}$ is the intersection of
$\Gamma_1$ and $\Gamma_2$.
\end{lemma}

\begin{proof}
For $\Gamma_1$, by setting $z=p_c+te^{i\pi/3}$, $0\le t\le
2(1-p_c)$,
\begin{equation}
\begin{split}
  F_1(t) &:= Re(f(p_c+te^{i\pi/3})) \\
  & = -\mu (p_c+\frac12 t -\qq)
+ \frac12 \ln\bigl(p_c^2+p_ct+t^2\bigr)
  - \frac1{2\gamma^2} \ln\bigl((1-p_c)^2-(1-p_c)t+t^2\bigr)
\end{split}
\end{equation}
and
\begin{equation}
  F_1'(t) = -\frac{t^2\bigl( (\gamma+1)^2t^2-(\gamma^2-1)t+2\gamma \bigr)}
  {2\gamma^2\bigl(p_c^2+p_ct+t^2\bigr)
  \bigl((1-p_c)^2-(1-p_c)t+t^2\bigr)}.
\end{equation}
The denominator is equal to $2\gamma^2|z|^2|1-z|^2$, and hence is
positive. To show that the numerator is positive, set
\begin{equation}
  T_1(t):= (\gamma+1)^2t^2-(\gamma^2-1)t+2\gamma.
\end{equation}
A simple calculus shows that
\begin{equation}
  \min_{ t\in [0, 2(1-p_c)]} T_1(t)= \begin{cases}
    T_1\bigl(\frac{\gamma^2-1}{2(\gamma+1)^2}\bigr), \qquad &1\le \gamma \le
    5, \\
    T_1(2(1-p_c)), \qquad &\gamma \ge 5.
  \end{cases}
\end{equation}
But $T_1\bigl(\frac{\gamma^2-1}{2(\gamma+1)^2}\bigr) \ge 2$
for $1\le \gamma\le 5$, and $T_1(2(1-p_c))=6$ for $\gamma\ge 5$,
and hence we find that $T_1(t)>0$ for $t\in[0, 2(1-p_c)]$ and for
all $\gamma\ge 1$. Thus we find that $F_1(t)$ is an increasing
function in $t\in [0, 2(1-p_c)]$.

For $\Gamma_2$, by setting $z= p_c+2(1-p_c)e^{i\pi/3} + x$, $x\ge
0$,
\begin{equation}
\begin{split}
  F_2(x) &:= Re\bigl(f (p_c+2(1-p_c)e^{i\pi/3} + x)\bigr) \\
  &= -\mu(1+x-q) + \frac12 \ln \bigl( (1+x)^2+3(1-p_c)^2\bigr)
  - \frac1{2\gamma^2} \ln\bigl( x^2+3(1-p_c)^2\bigr),
\end{split}
\end{equation}
and
\begin{equation}
  F'_2(x) = -\frac{T_2(x)}{\gamma^2(\gamma+1)^2 \bigl( (1+x)^2+3(1-p_c)^2\bigr)
  \bigl( x^2+3(1-p_c)^2\bigr)}
\end{equation}
where
\begin{equation}
\begin{split}
  T_2(x) =& (\gamma+1)^4x^4 +
  (\gamma^4+6\gamma^3+12\gamma^2+10\gamma+3)x^3 \\
  & + (2\gamma^3+13\gamma^2+20\gamma+9)x^2
  + 2(2\gamma^2+7\gamma+5)x + 6(\gamma+2)>0, \quad x\ge 0.
\end{split}
\end{equation}
Hence $F_2(x)$ is decreasing for $x\ge 0$.

For $\Gamma_3$, setting $z=1+R+iy$,  $0\le y\le
\sqrt{3}(1-p_c)$,
\begin{equation}
\begin{split}
  F_3(y) &= Re\bigl(f (1+R+iy)\bigr) \\
  &= -\frac{\mu\epsilon}{\nu M^{1/3}} -\mu(R+1-p_c) + \frac12 \ln\bigl( (1+R)^2+y^2\bigr)
  - \frac1{2\gamma^2} \ln\bigr( R^2+y^2\bigr).
\end{split}
\end{equation}
As $1\le \mu=(\gamma+1)^2/\gamma^2\le 4$ for $\gamma\ge 1$ and
$1-p_c=1/(\gamma+1)\le 1$, $F_3(y)$ can be made arbitrarily small
when $R$ is taken to be large. But
\begin{equation}
\begin{split}
   & Re (f(p_*))= Re \bigl(f(p_c+2(1-p_c)e^{i\pi/3})\bigr) \\
   &\quad = -\frac{\mu\epsilon}{\nu M^{1/3}} -\mu (1-p_c) + \frac12\ln\bigl( p_c^2+ 2p_c(1-p_c) +
   4(1-p_c)^2\bigr)
   - \frac1{2\gamma^2} \ln \bigl(5(1-p_c)^2\bigr)
\end{split}
\end{equation}
is bounded for $\gamma$ in a compact subset of $[1,\infty)$. Thus
the result follows.
\end{proof}

As $\gamma$ is in a compact subset of $[1,\infty)$, we assume from
here on that
\begin{equation}
   1\le \gamma\le \gamma_0
\end{equation}
for a fixed $\gamma_0\ge 1$. Now we split the contour
$\Gamma=\Gamma'\cup\Gamma''$ where $\Gamma'$ is the part of
$\Gamma$ in the disk $|z-p_c| < \delta$ for some $\delta>0$ which
will be chosen in the next paragraph, and $\Gamma''$ is the rest
of $\Gamma$. Let $\Gamma_{\infty}'$ be the image of the contour
$\Gamma'$ under the map $z\mapsto \nu M^{1/3}(z-p_c)$ and let
$\Gamma_{\infty}''=\Gamma_\infty\setminus \Gamma_\infty'$. Set
\begin{equation}
  \mathcal{H}(u)=\mathcal{H}'(u)+ \mathcal{H}''(u),
  \qquad \mathcal{H}_\infty(u)=\mathcal{H}_\infty'(u)+ \mathcal{H}_\infty''(u)
\end{equation}
where $\mathcal{H}'(u)$ (resp. $\mathcal{H}_\infty'(u)$) is the
part of the integral formula of $\mathcal{H}(u)$ (res.
$\mathcal{H}_\infty(u)$) integrated only over the contour $\Gamma'$
(resp. $\Gamma_{\infty}'$).

Fix $\delta$ such that
\begin{equation}\label{eq:delta}
  0< \delta < \min\bigl\{ \frac{\nu^3}{6C_0},
  \frac{1}{2(1+\gamma_0)} \bigr\},
  \qquad C_0:= 4^3+4(1+\gamma_0)^4.
\end{equation}
For $|s-p_c|\le \delta$,
\begin{equation}
\begin{split}
  \biggl|\frac1{4!}f^{(4)}(s) \biggr| = \frac14 \biggl|
  \frac{-1}{s^4}+\frac{\gamma^{-2}}{(1-s)^4} \biggr|
  &\le \frac14 \biggl(
  \frac1{(p_c-\delta)^4}+\frac{\gamma^{-2}}{(1-p_c-\delta)^4}
  \biggr) \\
  &\le \frac1{4}\bigl( 4^4+2^4(1+\gamma_0)^4 \bigr)=C_0.
\end{split}
\end{equation}
Hence by the Taylor's theorem, for $|z-p_c|\le \delta$,
\begin{equation}\label{eq:fineqmain}
\begin{split}
   \biggl|Re\bigl( f(z)-f(p_c) - \frac{f^{(3)}(p_c)}{3!}(z-p_c)^3 \bigr) \biggr|
  &\le \bigl| f(z)-f(p_c) - \frac{f^{(3)}(p_c)}{3!}(z-p_c)^3 \bigr|
  \\
  &   \le \biggl( \max_{|s-p_c|\le \delta}
\frac{|f^{(4)}(s)|}{4!} \biggr)
  |z-p_c|^4  \\
  &  \le C_0 |z-p_c|^4
  \le \frac{\nu^3}6 |z-p_c|^3.
\end{split}
\end{equation}
Therefore, by recalling \eqref{eq:fthreep}, we find for $0\le t\le
\delta$,
\begin{equation}
  Re\bigl( f(p_c+ t e^{i\pi/3}) \bigr) -f(p_c)
  \le - \frac{\nu^3}{6} t^3.
\end{equation}
Especially, from the Lemma \ref{lem:Refdecreasing} (note
that $Re(f(z))= Re(f(\overline{z}))$),
\begin{equation}\label{eq:finqcubic}
  \max_{z\in\Gamma''} Re(f(z)) \le Re\bigl(f(p_c+\delta e^{i\pi/3})\bigr)
  \le f(p_c) - \frac{\nu^3}{6}\delta^3.
\end{equation}

In the following sub-subsections, we consider two cases
separately; first when $u$ is in a compact subset of $\mathbb{R}$,
and the other when $u>0$.

\subsubsection{When $u$ is in a compact subset of $\mathbb{R}$.}

Suppose that $|u|\le u_0$ for some $u_0>0$. First we estimate
\begin{equation}\label{eq:124}
  |Z_M\mathcal{H}''(u)|\le \frac{|g(p_c)|}{2\pi (\nu M^{1/3})^{k-1}}
  \int_{\Gamma''}
  e^{-\nu M^{1/3} u Re(z-\qq)} e^{M Re(f(z)-f(p_c))} \frac1{|p_c-z|^k|g(z)|} |dz|.
\end{equation}
Using \eqref{eq:finqcubic} and $Re(z-\qq) \le \sqrt{(1+R)^2+3}$
for $z\in \Gamma''$,
\begin{equation}\label{eq:125}
  |Z_M \mathcal{H}''(u)| \le \frac{|g(p_c)|}{2\pi (\nu M^{1/3})^{k-1}}
    e^{\nu M^{1/3}u_0 \sqrt{(1+R_1)^2+3}}
    e^{-M\frac{\nu^3}6 \delta^3} \frac{L_{\Gamma}C_g}{\delta^k},
\end{equation}
where $L_\Gamma$ is the length of $\Gamma$ and $C_g>0$ is a
constant such that
\begin{equation}
  \frac1{C_g} \le \min_{z\in\Gamma''} |g(z)| \le C_g.
\end{equation}
For $\gamma \in [1,\gamma_0]$, $L_\Gamma$, $C_g$ and $|g(p_c)|$
are uniformly bounded, and hence
\begin{equation}\label{eq:estucomp1}
  | Z_M \mathcal{H}''(u)| \le e^{-\frac{\nu^3}{12} \delta^3 M}
\end{equation}
when $M$ is sufficiently large.

On the other hand, for $a\in \Gamma_\infty''$,
we have $a=t e^{\pm i\pi/3}$,
$\delta\nu M^{1/3} \le t <+\infty$, and hence
\begin{equation}\label{eq:estucomp2}
\begin{split}
  | \mathcal{H}_{\infty}''(u) |
  & = \biggl| \frac{e^{-\epsilon u}}{2\pi} \int_{\Gamma_\infty''}
  \frac{e^{-ua} e^{\frac13 a^3}}{a^k} da \biggr|
  \le \frac{e^{\epsilon u_0}}{2\pi} \int_{\Gamma_{\infty}''}
  \frac{e^{u_0 |a|} e^{\frac13 Re(a^3)}}{|a|^k} |da|  \\
  &\le \frac{e^{\epsilon u_0}}{\pi} \int_{\delta\nu M^{1/3}}^\infty
  \frac{e^{u_0t-\frac13 t^3} }{t^k} dt
  \le e^{-\frac{\nu^3}{6}\delta^3M}
\end{split}
\end{equation}
when $M$ is sufficiently large.

Now we estimate
$| Z_M \mathcal{H}'(u) - e^{-\epsilon u}\mathcal{H}_\infty'(u) |$. Using the change of
variables $a=\nu M^{1/3}(z-p_c)$ for the integral \eqref{eq:Hinf1} for $\mathcal{H}_\infty'(u)$,
\begin{equation}\label{eq:Htemppart}
\begin{split}
  &| Z_M \mathcal{H}'(u) - \mathcal{H}_\infty'(u) | \\
&\quad \le \frac{\nu M^{1/3}}{2\pi} \int_{\Gamma'} \frac{e^{-\nu
M^{1/3} uRe(z-\qq)}}{|\nu M^{1/3}(z-p_c)|^k} \biggl|
e^{M(f(z)-f(p_c)} \frac{g(p_c)}{g(z)} - e^{M\frac{\nu^3}3
(z-p_c)^3} \biggr| |dz|.
\end{split}
\end{equation}
We split the integral into two parts. Let
$\Gamma'=\Gamma'_1\cup\Gamma_2'$ where
$\Gamma_1'=\Gamma_0\cup\overline{\Gamma_0}$ and
$\Gamma_2'=\Gamma'\setminus\Gamma_1'$.

For $z\in\Gamma_1'$, $|z-p_c|=\epsilon /(2\nu M^{1/3})$.
Hence using \eqref{eq:fineqmain},
\begin{equation}\label{eq:estGamma1p}
\begin{split}
  &|e^{M(f(z)-f(p_c))} - e^{M\frac{\nu^3}3(z-p_c)^3} | \\
  &\quad  \le \max \bigl( |e^{M(f(z)-f(p_c))}|, |e^{M\frac{\nu^3}3(z-p_c)^3}| \bigr)
  \cdot M|f(z)-f(p_c)-\frac{\nu^3}3 (z-p_c)^3 | \\
  &\quad \le e^{MRe(\frac{\nu^3}3(z-p_c)^3+\frac{\nu^3}6|z-p_c|^3 )} MC_0 |z-p_c|^4 \\
  &\quad \le e^{\frac1{16}\epsilon^3} \frac{C_{0}\epsilon^4}{16 \nu^4 M^{1/3}}.
\end{split}
\end{equation}
Also
\begin{equation}
  \biggl| \frac{g(p_c)}{g(z)}-1 \biggr|
  \le \frac1{|g(z)|} \displaystyle \max_{|s-p_c|\le
\frac{\epsilon}{2\nu M^{1/3}}} |g'(s)|
\cdot |z-p_c|
  \le \frac{C_0\overline{C}\epsilon}{2\nu M^{1/3}}
\end{equation}
where
\begin{equation}
  \overline{C}:= \max\{
|g'(s)| : |s-p_c|\le
\frac{\epsilon}{2\nu M^{1/3}}, s\in\Gamma_2'\}
\end{equation}
which is uniformly bounded as $\Gamma$ is uniformly away
from the singular points of $g$.
Hence using
\begin{equation}
\begin{split}
  &\biggl| e^{M(f(z)-f(p_c)} \frac{g(p_c)}{g(z)} - e^{M\frac{\nu^3}3 (z-p_c)^3} \biggr| \\
  &\quad =  \biggl| \bigl( e^{M(f(z)-f(p_c)} - e^{M\frac{\nu^3}3(z-p_c)^3}\bigr)
\frac{g(p_c)}{g(z)}
  +e^{M\frac{\nu}3 (z-p_c)^3} \bigl( \frac{g(p_c)}{g(z)}-1 \bigr)\biggr|  \\
  &\quad \le \biggl(
  \frac{C_0\epsilon^4}{16\nu^4}e^{\frac1{16} \epsilon^3} +
\frac{C_0\overline{C}}{2\nu}
\biggr)\cdot \frac1{M^{1/3}},
\end{split}
\end{equation}
and the fact that the length of $\Gamma_1'$ is $4\pi R_0/(3M^{1/3})$,
we find that the part of the integral in \eqref{eq:Htemppart}
over $\Gamma_1'$ is less than or equal to some constant divided by $M^{1/3}$.

For $z\in \Gamma_2'$, we have
$z=p_c + te^{\pm i\pi/3}$, $\epsilon/(2\nu M^{1/3}) \le t\le \delta$.
From \eqref{eq:fineqmain} (cf. \eqref{eq:estGamma1p})
\begin{equation}
\begin{split}
  &|e^{M(f(z)-f(p_c))} - e^{M\frac{\nu^3}3(z-p_c)^3} | \\
  &\quad \le e^{MRe(\frac{\nu^3}3(z-p_c)^3+\frac{\nu^3}6|z-p_c|^3 )}
MC_{0} |z-p_c|^4 \\
  &\quad \le e^{-M\frac{\nu^3}6 t^3} \cdot C_{0} Mt^4.
\end{split}
\end{equation}
Also
\begin{equation}
  \biggl|  \frac{g(p_c)}{g(z)}-1 \biggr|
  \le \frac{1}{|g(z)|} \max_{s\in \Gamma_2'} |g'(s)| \cdot |p_c-z|
  \le C_0\overline{C} t,
\end{equation}
and hence
\begin{equation}
\begin{split}
  &\biggl| e^{M(f(z)-f(p_c)} \frac{g(p_c)}{g(z)} -
e^{M\frac{\nu^3}3 (z-p_c)^3} \biggr| \\
  &\quad =  \biggl| \bigl( e^{M(f(z)-f(p_c)} - e^{M\frac{\nu^3}3(z-p_c)^3}\bigr)
\frac{g(p_c)}{g(z)} \biggr|
  +
\biggl| e^{M\frac{\nu^3}3 (z-p_c)^3} \bigl( \frac{g(p_c)}{g(z)}-1 \bigr)\biggr|  \\
  &\quad
\le e^{-\frac{\nu^3}6 Mt^3}C_0^3Mt^4 + e^{-\frac{\nu^3}3 Mt^3} C_0\overline{C}t
\le (C_0^3+C_0\overline{C}) e^{-\frac{\nu^3}6 Mt^3} (Mt^4+t).
\end{split}
\end{equation}
Using
\begin{equation}
  e^{-\nu M^{1/3} uRe (z-\qq)}  \le e^{\nu M^{1/3}u_0 (|p_c-\qq|+|z-p_\qq|)}
  = e^{\epsilon u_0+\nu u_0 M^{1/3}t},
\end{equation}
we find by substituting $z=p_c+ t e^{\pm i\pi/3}$ into \eqref{eq:Htemppart},
the part of the integral in \eqref{eq:Htemppart} over $\Gamma_2'$ is less than or equal to
\begin{equation}
  \frac{\nu M^{1/3}}{\pi}(C_0^3+C_0\overline{C})
\int_{\frac{\epsilon}{2\nu M^{1/3}}}^{\delta}
  \frac{e^{\epsilon u_0+\nu u_0 M^{1/3}t}}{(\nu M^{1/3} t)^k}
e^{-\frac{\nu^3}6Mt^3} (Mt^4+t) dt.
\end{equation}
Then by the change of variables $s=\nu M^{1/3}t$, the last integral is less than or equal to
\begin{equation}\label{eq:137}
  \frac{(C_0^3+C_0\overline{C})e^{\epsilon u_0}}{\pi M^{1/3}}
  \int_{\epsilon/2}^\infty
   \frac{e^{u_0s-\frac16 s^3}}{s^k} \bigl( \frac1{\nu^4}s^4+\frac1{\nu}s \bigr)
ds,
\end{equation}
which is a constant divided by $M^{1/3}$.
Therefore, together with the estimate for the part over $\Gamma_1'$,
this implies that
\begin{equation}\label{eq:estucomp3}
   | Z_M \mathcal{H}'(u) - \mathcal{H}_\infty'(u) |
  \le \frac{C_1}{M^{1/3}}
\end{equation}
for some positive constant $C_1>0$ for any $M>0$. Now combining
\eqref{eq:estucomp1}, \eqref{eq:estucomp2} and
\eqref{eq:estucomp3}, we find that for any $u_0>0$, there are
constants $C>0, M_0>0$ which may depend on $u_0$ such that
\begin{equation}\label{eq:Hdiffcompact}
   | Z_M \mathcal{H}(u) - \mathcal{H}_\infty(u) |
  \le \frac{C}{M^{1/3}}
\end{equation}
for $|u|\le u_0$ and $M\ge M_0$.
obtain \eqref{eq:estprop1}.

\subsubsection{When $u>0$.}

We first estimate $|Z_M \mathcal{H}''(u)|$ using \eqref{eq:124}.
For $z\in \Gamma''$, $Re(z-\qq) = Re(p_c-\qq)+Re(z-p_c) \ge
\frac{\epsilon}{\nu M^{1/3}} + \frac12\delta$, and hence as $u>0$,
\begin{equation}\label{eq:140}
  |Z_M \mathcal{H}''(u)| \le \frac{|g(p_c)|}{2\pi (\nu M^{1/3})^{k-1}}
    e^{-\epsilon u} e^{-\frac12 \nu M^{1/3}\delta u}
    e^{-M\frac{\nu^3}6 \delta^3} \frac{L_{\Gamma}}{\delta^k C_g}.
\end{equation}

On the other hand, for $a\in \Gamma_\infty''$, by estimating as in
\eqref{eq:estucomp2} but now using $u>0$ and $Re(a)\ge
\frac12\delta \nu M^{1/3}$ for $a\in\Gamma_\infty''$, we find
\begin{equation}\label{eq:141}
\begin{split}
  |\mathcal{H}_\infty''(u)|
  &\le \frac{e^{-\epsilon u}}{2\pi} \int_{\Gamma_\infty''}
  \frac{e^{- Re(a)+\frac13 Re(a^3)}}{|a|^k} |da|
  \le \frac{e^{-\epsilon u}e^{-\frac12 \delta \nu M^{1/3} u}}{2\pi} \int_{\Gamma_\infty''}
  \frac{e^{\frac13 Re(a^3)}}{|a|^k} |da| \\
  & \le e^{-\epsilon u} e^{-\frac12 \delta \nu M^{1/3} u}
  e^{-M \frac{\nu^3}{6} \delta^3}
\end{split}
\end{equation}
when $M$ is sufficiently large.

In order to estimate $| Z_M \mathcal{H}'(u) -
\mathcal{H}_\infty'(u) |$, we note that as $u>0$, for $z\in
\Gamma'$,
\begin{equation}
   e^{-\nu M^{1/3} u Re(z-\qq)}
  \le e^{-\frac12 \epsilon u}
\end{equation}
which is achieved at $z=p_c-\frac{\epsilon}{2\nu M^{1/3}}$.
Using the same estimates as in the case when $u$ is in a compact set for
the rest of the terms of \eqref{eq:Htemppart}, we find that
\begin{equation}\label{eq:143}
\begin{split}
  | Z_M \mathcal{H}'(u) - \mathcal{H}_\infty'(u) | \le
  \frac{C_2e^{-\frac12 \epsilon u}}{M^{1/3}}
\end{split}
\end{equation}
for some constant $C_2>0$ when $M$ is large enough.
Now \eqref{eq:140}, \eqref{eq:141} and \eqref{eq:143} yield
\eqref{eq:estprop1} for the case when $u>0$.

\subsection{Proof of Proposition \ref{prop:HJmaina} (ii)}\label{sec:Jmain}

\begin{figure}[ht]
\centerline{\epsfxsize=11cm\epsfbox{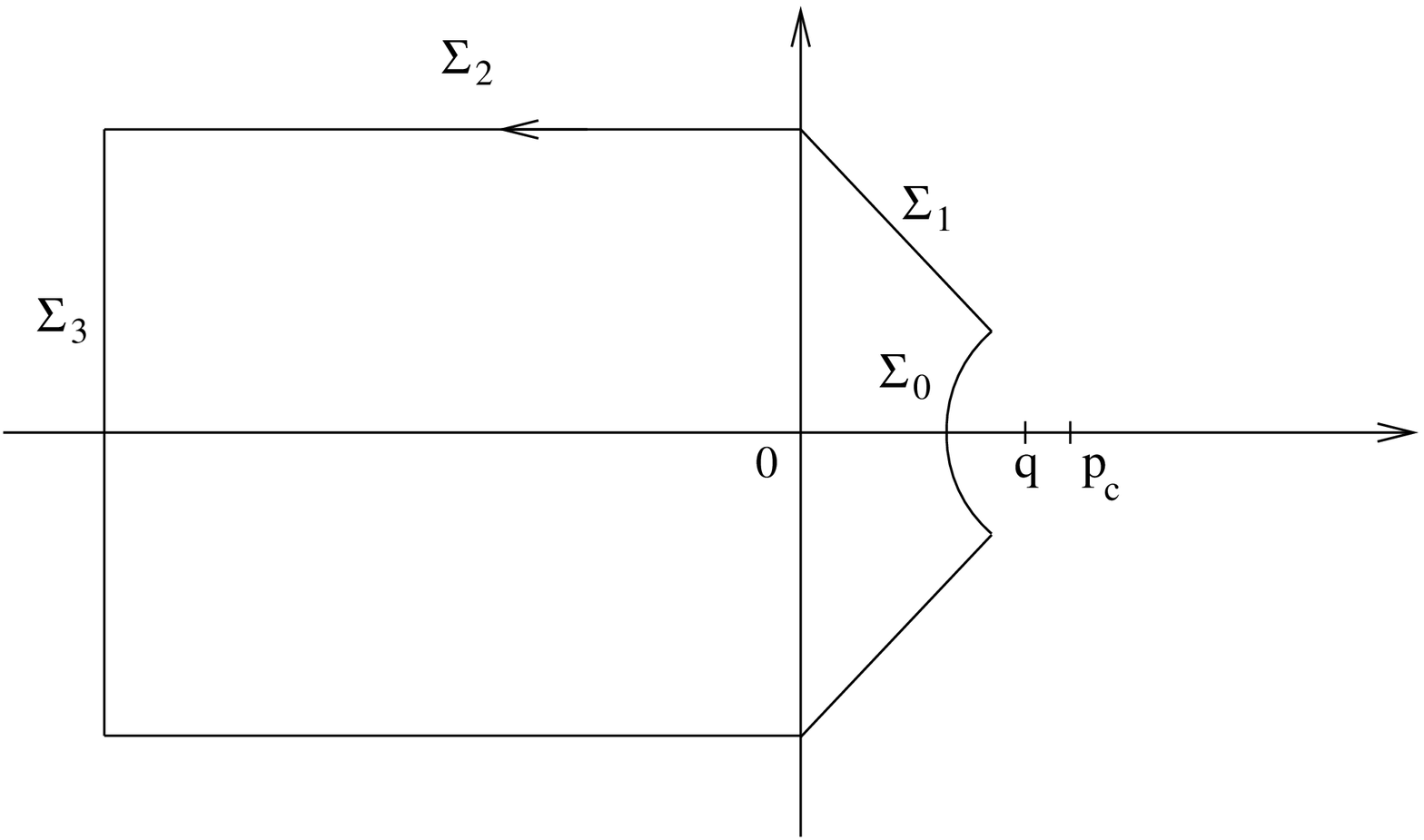}} \caption{Contour
$\Sigma$}\label{fig:Steep2}
\end{figure}

We choose the contour $\Sigma$ explicitly, which will be a
steep-descent curve of $-f(z)$, as follows. Let $R>0$. Define
\begin{eqnarray}
  \Sigma_0 &:=& \{ p_c+ \frac{3\epsilon}{\nu M^{1/3}}e^{i(\pi-\theta)} :
0\le \theta\le  \frac\pi3\} \\
  \Sigma_1 &:=& \{ p_c+ te^{2i\pi/3} :
\frac{3\epsilon}{\nu M^{1/3}} \le t \le 2 p_c
  \} \\
  \Sigma_2 &:=& \{ p_c+2 p_ce^{2i\pi/3} - x : 0\le x\le R\} \\
  \Sigma_3 &:=& \{ -R+i(\sqrt{3} p_c- y) : 0\le y\le \sqrt{3}p_c \},
\end{eqnarray}
and set
\begin{equation}
  \Sigma= \bigl( \cup_{k=0}^3 \Sigma_k \bigr)
\cup \overline{\bigl( \cup_{k=0}^3 \Sigma_k \bigr)}.
\end{equation}
The orientation of $\Sigma$ is counter-clockwise.
See Figure \ref{fig:Steep2}.
We first prove decay properties of $Re(-f(z))$ analogous
to Lemma \ref{lem:Refdecreasing}.

\begin{lemma}\label{lem:Refminus}
  For $\gamma\ge 1$, $Re(-f(z))$ is decreasing for
$z\in\Sigma_1\cup\Sigma_2$ as $Re(z)$ decreases. Also for $\gamma$
in a compact subset of $[1,\infty)$, we can take large $R>0$
(independent of $\gamma$) such that
\begin{equation}
  \max_{z\in\sigma_3} Re(-f(z)) \le Re
\bigl(-f(p_*))\bigr).
\end{equation}
where $p_*=p_c+2p_ce^{2i\pi/3}$ is the intersection
of $\Sigma_1$ and $\Sigma_2$.
\end{lemma}

\begin{proof}
For $z\in \Sigma_1$, by setting $z=p_c+t e^{2i\pi/3}$, $0\le t\le 2p_c$,
\begin{equation}
\begin{split}
  F_1(t)& := Re (-f(p_c+t e^{2i\pi/3}))  \\
  & = \mu(p_c-\frac12 t -c) - \frac12 \ln (p_c^2-p_c t+t^2)
  + \frac1{2\gamma^2} \ln ((1-p_c)^2+(1-p_c)t +t^2)
\end{split}
\end{equation}
and hence
\begin{equation}
  F_1'(t) = -\frac{t^2\bigl( (\gamma+1)^2t^2+(\gamma^2-1)t+2\gamma \bigr)}
  {2\gamma^2\bigl(p_c^2-p_ct+t^2\bigr)
  \bigl((1-p_c)^2+(1-p_c)t+t^2\bigr)},
\end{equation}
which is non-negative for all $t \ge 0$.

For $\Sigma_2$, by setting $z= p_c+2p_c e^{2i\pi/3} - x$, $x\ge 0$,
\begin{equation}
\begin{split}
  F_2(x) &:= Re\bigl(- f (p_c+2p_c e^{2i\pi/3} - x)\bigr) \\
  &= \mu(-x-c) - \frac12 \ln \bigl( x^2+3p_c^2\bigr)
  + \frac1{2\gamma^2} \ln\bigl( (1+x)^2+3p_c^2\bigr).
\end{split}
\end{equation}
A direct computation shows that
\begin{equation}
  F'_2(x) = -\frac{T_2(x)}{\gamma^2(\gamma+1)^2 \bigl( x^2+3p_c^2\bigr)
  \bigl( (1+x)^2+3p_c^2\bigr)}
\end{equation}
where
\begin{equation}
\begin{split}
  T_2(x) =& (\gamma+1)^4x^4 +
  (3\gamma^4+10\gamma^3+12\gamma^2+6\gamma+1)x^3 \\
  & + (9\gamma^4+20\gamma^3+13\gamma^2+2\gamma)x^2
  + 2(5\gamma^4+7\gamma^3+2\gamma^2)x + 6(2\gamma^4+\gamma^3).
\end{split}
\end{equation}
Hence $F_2(x)$ is decreasing for $x\ge 0$.

For $\Sigma_3$, setting $z=-R+i(\sqrt{3}p_c-y)$,  $0\le y\le
\sqrt{3}p_c$,
\begin{equation}
\begin{split}
  F_3(y) &= Re\bigl(-f (-R+i(\sqrt{3} p_c -y)\bigr) \\
  &=  \mu(-R-c) - \frac12 \ln\bigl( R^2+(\sqrt{3}p_c-y)^2\bigr)
  + \frac1{2\gamma^2} \ln\bigr( (1+R)^2+(\sqrt{3} p_c-y)^2\bigr).
\end{split}
\end{equation}
When $R\to +\infty$, $F_3(y)$ can be made arbitrarily small. But
\begin{equation}
\begin{split}
   & Re \bigl(-f(p_c+2p_c e^{2i\pi/3})\bigr) \\
   &\quad =  -\mu c - \ln\bigl( \sqrt{3} p_c \bigr)
   + \frac1{2\gamma^2} \ln \bigl(1+3 p_c^2\bigr)
\end{split}
\end{equation}
is bounded for $\gamma$ in a compact subset of $[1,\infty)$. Thus
the result follows.
\end{proof}

Let $\delta$ be given in \eqref{eq:delta}.
Let $\Sigma=\Sigma'\cup\Sigma''$ where $\Sigma'$ is the part of $\Sigma$
that lies in the disk $|z-p_c|<\delta$,
and let $\Sigma''=\Sigma\setminus \Sigma'$. Let $\Sigma_\infty'$ be the image
of $\Sigma'$ under the map $z\mapsto \nu M^{1/3}(z-p_c)$
and let $\Sigma_\infty'' = \Sigma_\infty\setminus \Sigma_\infty'$. Set
\begin{equation}
  \mathcal{J}(v)=\mathcal{J}'(v) + \mathcal{J}''(v),
\qquad \mathcal{J}_\infty=\mathcal{J}_\infty'(v) + \mathcal{J}_\infty''(v)
\end{equation}
where $\mathcal{J}'(v)$ (resp. $\mathcal{J}_\infty'(v)$) is the part of
the integral formula of $\mathcal{J}(v)$ (resp. $\mathcal{J}_\infty(v)$)
integrated over the contour $\Sigma'$ (resp. $\Sigma_\infty'$).

As before, we consider two cases separately; first case when $v$ is in
a compact subset of $\mathbb{R}$, and the second case when $v>0$.

\subsubsection{When $v$ is in a compact subset of $\mathbb{R}$.}

There is $v_0>0$ such that $|v|\le v_0$. First, we estimate
\begin{equation}
  \bigl| \frac1{Z_M} \mathcal{J}''(v) \bigr|
  \le \frac{(\nu M^{1/3})^{k+1}}{2\pi |g(p_c)|}
   \int_{\Sigma''} e^{\nu M^{1/3} v_0 |z-\qq|} e^{M Re(-f(z)+f(p_c))}
  |p_c-z|^k |g(z)| |dz|.
\end{equation}
From Lemma \ref{lem:Refminus} and \eqref{eq:fineqmain},
following the proof of \eqref{eq:finqcubic},
we find
\begin{equation}
  \max_{z\in \Sigma''} Re(-f(z))
  \le Re(-f(p_c+\delta e^{2i\pi/3})
  \le f(p_c)- \frac{\nu^3}6 \delta^3.
\end{equation}
Hence using $|z-\qq|\le \sqrt{(R_1+1)^2+3}$ for $z\in \Sigma''$,
\begin{equation}
  \bigl| \frac1{Z_M} \mathcal{J}''(v) \bigr|
  \le \frac{(\nu M^{1/3})^{k+1}}{2\pi |g(p_c)|}
  e^{\nu M^{1/3} v_0\sqrt{(R_1+1)^2+3}} e^{-\frac{\nu^3}6\delta^3M}
  \delta^k \widetilde{C}_g L_\Sigma,
\end{equation}
where $\widetilde{C}_g$ is the maximum of $|g(z)|$ over $z\in\Sigma''$
and $L_\Sigma$ is the length of $\Sigma''$, both of which are
uniformly bounded. Hence we find that when $M$ is sufficiently large,
\begin{equation}\label{eq:161}
  \bigl| \frac1{Z_M} \mathcal{J}''(v) \bigr|
  \le e^{-\frac{\nu^3}{12} \delta^3 M}.
\end{equation}

When $a\in \Sigma_\infty''$, $a= te^{\pm 2i\pi/3}$,
$\delta \nu M^{1/3} \le t<+\infty$, and
\begin{equation}\label{eq:162}
   | \mathcal{J}_\infty''(v)|
   = \biggl| \frac{e^{\epsilon v}}{2\pi}
\int_{\Sigma_\infty''} e^{va-\frac13a^3}a^k da   \biggr|
   \le \frac{e^{\epsilon v_0}}{\pi}  \int_{\delta \nu M^{1/3}}^\infty
e^{v_0t-\frac13 t^3} t^k dt
   \le e^{-\frac{\nu^3}{6}\delta^3 M}
\end{equation}
when $M$ is sufficiently large.

Finally we estimate
\begin{equation}
\begin{split}
  & \bigl| \frac1{Z_M} \mathcal{J}'(v) -\mathcal{J}_\infty'(v) \bigr| \\
&
\le \frac{\nu M^{1/3}}{2\pi} \int_{\Sigma'}
|e^{\nu M^{1/3} v(z-\qq)} | |\nu M^{1/3}(z- p_c)|^k
\biggl| e^{-M(f(z)-f(p_c))} \frac{g(z)}{g(p_c)} - e^{-M\frac{\nu^3}3(z-p_c)^3}
\biggr| |dz|.
\end{split}
\end{equation}
As before, we split the contour $\Sigma'=\Sigma_1'\cup\Sigma_2'$ where
$\Sigma_1'=\Sigma_0\cup\overline{\Sigma_0}$ and
$\Sigma_2'=\Sigma'\setminus\Sigma_1'$, and by following the steps
of \eqref{eq:estGamma1p}-\eqref{eq:137}, we arrive at
\begin{equation}\label{eq:164}
\begin{split}
  \bigl| \frac1{Z_M} \mathcal{J}'(v) -\mathcal{J}_\infty'(v) \bigr|
\le \frac{C_3}{M^{1/3}},
\end{split}
\end{equation}
for some constant $C_3>0$ when $M$ is large enough. From
\eqref{eq:161}, \eqref{eq:162} and \eqref{eq:164}, we obtain
\eqref{eq:estprop3}.

\subsubsection{When $v>0$.}

The proof in this case is again very similar to the
estimate of $\mathcal{H}(u)$ when $u>0$. Then only
change is the following estimates:
\begin{equation}
  Re(z-\qq) \le Re(p_c+\delta e^{2i\pi/3}-\qq)
= \frac{\epsilon}{\nu M^{1/3}} - \frac12 \delta,
\qquad z\in \Sigma'',
\end{equation}
\begin{equation}
  Re(z-\qq)\le Re(p_c+ \frac{3\epsilon}{\nu M^{1/3}} e^{2i\pi/3}-\qq)
  = -\frac{\epsilon}{2\nu M^{1/3}},
\qquad z\in \Sigma',
\end{equation}
and
\begin{equation}
  Re(z-\qq) = -\frac12 |z-p_c| + \frac{\epsilon}{\nu M^{1/3}},
\qquad z\in \Sigma_1.
\end{equation}

Then for large enough $M>0$,
\begin{equation}
   \biggl| \frac1{Z_M} \mathcal{J}''(v) \biggr|
\le e^{\epsilon v} e^{-\frac12 \delta \nu M^{1/3} v}
e^{-\frac{\nu^3}{12} \delta^3 M},
\end{equation}
\begin{equation}
   \bigl| \mathcal{J}_\infty''(v) \bigr|
\le e^{\epsilon v} e^{-\frac12 \delta \nu M^{1/3} v}
e^{-\frac{\nu^3}{6} \delta^3 M},
\end{equation}
and
\begin{equation}
   \biggl| \frac1{Z_M} \mathcal{J}'(v)- \mathcal{J}_\infty'(v) \biggr|
\le \frac{C}{M^{1/3}} e^{-\frac12 \epsilon v}
\end{equation}
for some constant $C>0$.
We skip the detail.

\subsection{Proof of Theorem \ref{thm:main} (a)}

From the Proposition \ref{prop:HJmaina}, the discussion on the
subsection \ref{sec:basicas} implies that under the assumption of
Theorem \ref{thm:main} (a),
\begin{equation}\label{eq:Proofatemp}
  \mathbb{P} \biggl( \bigr( \lambda_1 - (1+\gamma^{-1})^2 \bigr) \cdot
  \frac{\gamma}{(1+\gamma)^{4/3}} M^{2/3} \le x \biggr)
\end{equation}
converges, as $M\to\infty$, to the Fredholm determinant of the
operator acting on $L^2((0,\infty))$ whose kernel is
\begin{equation}
  \int_0^\infty
  \mathcal{H}_\infty(x+u+y)\mathcal{J}_\infty(x+v+y)dy.
\end{equation}
From the integral representation \eqref{eq:Airyintrep} of the Airy
function, by simple changes of variables,
\begin{equation}
  Ai(u) = \frac{-1}{2\pi i} \int_{\Gamma_\infty} e^{-ua +\frac13 a^3} da
= \frac1{2\pi i} \int_{\Sigma_\infty} e^{ub - \frac13 b^3} db,
\end{equation}
and hence a simple algebra shows that
\begin{equation}\label{eq:KinAste}
\begin{split}
  &  \int_0^\infty
  \mathcal{H}_\infty(u+y)\mathcal{J}_\infty(v+y)dy-
  \int_0^\infty
  e^{-\epsilon (u+y)}Ai(u+y) Ai(v+y) e^{\epsilon (v+y)} dy \\
  &\quad = \frac{e^{-\epsilon (u-v)}}{(2\pi)^2} \int_0^\infty dy
  \int_{\Gamma_\infty} da
  \int_{\Sigma_\infty} db \,
  e^{-(u+y)a+\frac13a^3} e^{(v+y)b-\frac13b^3}
  \biggl(\biggl(\frac{b}{a}\biggr)^k -1\biggr) \\
  &\quad
  = \sum_{m=1}^{k} \frac{e^{-\epsilon (u-v)}}{(2\pi)^2}
  \int_{\Gamma_\infty} da
  \int_{\Sigma_\infty} db \,
  e^{-ua+\frac13a^3} e^{vb-\frac13b^3}
  \frac{(b-a)b^{m-1}}{a^{m}} \int_0^\infty
  e^{-(a-b)y} dy \\
  &\quad
  = -\sum_{m=1}^{k} \frac{e^{-\epsilon (u-v)}}{(2\pi)^2}
  \int_{\Gamma_\infty} da
  \int_{\Sigma_\infty} db \,
  e^{-ua+\frac13a^3} e^{vb-\frac13b^3} \frac{b^{m-1}}{a^m} \\
  &\quad
  = \sum_{m=1}^{k} e^{-\epsilon u}s^{(m)}(u)t^{(m)}(v) e^{\epsilon v},
\end{split}
\end{equation}
where the choice of the contours $\Sigma_\infty$ and $\Gamma_\infty$
ensures that $Re(a-b)>0$ which is used in the third equality.

Let $\mathbf{E}$ be the multiplication operator by $e^{-\epsilon
u}$; $(\mathbf{E}f)(u) = e^{-\epsilon u} f(u)$. The computation
\eqref{eq:KinAste} implies that \eqref{eq:Proofatemp} converges to
\begin{equation}
  \det \bigl(1- \mathbf{E}\mathbf{A}_x \mathbf{E}^{-1} -\sum_{m=1}^k
  \mathbf{E} s^{(m)}_x \otimes t^{(m)}_x \mathbf{E}^{-1} \bigr).
\end{equation}
The general formula of the Fredholm determinant of a finite-rank
perturbation of an operator yields that this is equal to
\begin{equation}
  \det \bigl(1- \mathbf{E}\mathbf{A}_x \mathbf{E}^{-1} \bigr)
  \cdot \det \biggl( \delta_{mn} - < \frac{1}{1-\mathbf{E}\mathbf{A}_x
  \mathbf{E}^{-1}}  \mathbf{E} s^{(m)}_x , t^{(n)}_x \mathbf{E}^{-1} >
   \biggr)_{1\le m,n\le k},
\end{equation}
which is equal to \eqref{eq:Fkdef} due to the proof of the
following Lemma. This completes the proof.

\begin{lemma}\label{lem:Fk}
The function $F_k(x)$ in Definition \ref{def:Fk} is well-defined.
Also $s^{(m)}(u)$ defined in \eqref{eq:sdef} can be written as
\begin{equation}
  s^{(m)}(u) =  \sum_{\ell+3n = m-1} \frac{(-1)^n}{3^n \ell ! n!}
  u^\ell + \frac1{(m-1)!} \int_\infty^u (u-y)^{m-1} Ai(y)dy.
\end{equation}
\end{lemma}

\begin{proof}
It is known that $\mathbf{A}_x$ has norm less than $1$ and is
trace class (see, e.g. \cite{TracyWidom}). The only thing we need
to check is that the product $< \frac{1}{1-\mathbf{A}_x}
  s^{(m)}, t^{(n)} >$ is finite. By using the standard
steepest-descent analysis,
\begin{equation}\label{eq:Airyas}
  Ai(u) \sim \frac{1}{2\sqrt\pi u^{1/4}} e^{-\frac23
  u^{3/2}},  \qquad u\to +\infty.
\end{equation}
and
\begin{equation}\label{eq:tas}
  t^{(m)}(v) \sim \frac{v^{m/2}}{2\sqrt\pi v^{3/4}} e^{-\frac23
  v^{3/2}},  \qquad v\to +\infty.
\end{equation}
But for $s^{(m)}$, since the critical point $a=i$ is above the
pole $a=0$, the residue at $a=0$ contributes to the asymptotics
and $s^{(m)}(u)$ grows in powers of $u$ as $u\to+\infty$:
\begin{equation}\label{eq:sas}
  s^{(m)}(u)\sim  \sum_{\ell+3n = m-1} \frac{(-1)^n}{3^n \ell ! n!}
  u^\ell + \frac{(-1)^m}{2\sqrt\pi u^{m/2} u^{1/4}} e^{-\frac23
  u^{3/2}} , \qquad u\to +\infty.
\end{equation}
But the asymptotics \eqref{eq:Airyas} of the Airy function as
$u\to\infty$ implies that for any $U, V\in \mathbb{R}$, there is a
constant $C>0$ such that
\begin{equation}
  |\mathbf{A}(u,v) | \le C e^{-\frac23(u^{3/2}+v^{3/2})}, \qquad u\ge U,
  \quad v\ge V,
\end{equation}
which, together with \eqref{eq:tas}, implies that the inner
product $< \frac{1}{1-\mathbf{A}_x} s^{(m)}, t^{(n)} >$ is finite.

Also $s^{(m)}(u)$ defined in \eqref{eq:sdef} satisfies
$\frac{d^m}{du^m} s^{(m)}(u) = Ai(u)$. Hence $s^{(m)}(u)$ is
$m$-folds integral of $Ai(u)$ from $\infty$ to $u$ plus a
polynomial of degree $m-1$. But the asymptotics \eqref{eq:sas}
determines the polynomial and we obtain the result.
\end{proof}

\subsection{Proof of Theorem \ref{thm:main3}}
\label{sec:main3}

The analysis is almost identical to that of Proof of Theorem
\ref{thm:main} (a) with the only change of the scaling
\begin{equation}
  \pi_j^{-1}= 1+\gamma^{-1}-\frac{w_j}{M^{1/3}}.
\end{equation}
We skip the detail.

\section{Proofs of Theorem \ref{thm:main} (b)}
\label{sec:Gausslimit}

We assume that for some $1\le k\le r$,
\begin{equation}\label{eq:picon2}
  \pi^{-1}_1=\cdots = \pi^{-1}_k > 1+\gamma^{-1}
\end{equation}
are in a compact subset of $(1+\gamma^{-1}, \infty)$, and
$\pi^{-1}_{k+1}, \dots, \pi^{-1}_r$ are in a compact subset of
$(0, \pi_1^{-1})$.

For the scaling \eqref{eq:scalegeneral}, we take
\begin{equation}
  \alpha=1/2
\end{equation}
and
\begin{equation}\label{eq:mu2}
  \mu=\mu(\gamma) := \frac1{\pi_1}+\frac{\gamma^{-2}}{(1-\pi_1)}, \qquad
  \nu=\nu(\gamma) := \sqrt{ \frac1{\pi_1^2} - \frac{\gamma^{-2}}{(1-\pi_1)^2} }
\end{equation}
so that
\begin{equation}
  \xi= \mu+\frac{\nu x}{\sqrt{M}}.
\end{equation}
It is direct to check that the term inside the square-root of
$\nu$ is positive from the condition \eqref{eq:picon2}. Again, the
reason for such a choice will be clear during the subsequent
asymptotic analysis.

The functions \eqref{eq:Hcal} and \eqref{eq:Jcal} are now
\begin{equation}\label{eq:H2}
  \mathcal{H}(u)
  = \frac{\nu M^{1/2}}{2\pi}
  \int_\Gamma e^{-\nu M^{1/2} u(z-\qq)} e^{Mf(z)}
  \frac1{\bigl(\pi_1-z\bigr)^kg(z)} dz
\end{equation}
and
\begin{equation}\label{eq:J2}
  \mathcal{J}(v) = \frac{\nu M^{1/2}}{2\pi} \int_\Sigma
  e^{\nu M^{1/2} v(z-\qq)} e^{-Mf(z)}
  \bigl(\pi_1-z\bigr)^kg(z) dz
\end{equation}
where
\begin{equation}
   f(z) := -\mu (z-\qq) +\log(z) - \frac1{\gamma^2}\log(1-z),
\end{equation}
where log is the principal branch of logarithm, and
\begin{equation}
  g(z) := \frac{1}{(1-z)^{r}}
  \prod_{\ell=k+1}^r (\pi_\ell-z).
\end{equation}
The arbitrary parameter $\qq$ will be chosen in \eqref{eq:qq2}
below.
 Now as
\begin{equation}
  f'(z) = -\mu + \frac1{z} - \frac{1}{\gamma^2(z-1)},
\end{equation}
with the choice \eqref{eq:mu2} of $\mu$, two critical points of
$f$ are $z=\pi_1$ and $z=\frac1{\mu \pi_1}$. From the condition
\eqref{eq:picon2}, it is direct to check that
\begin{equation}
  \pi_1<\frac{\gamma}{1+\gamma}< \frac1{\mu\pi_1} <1.
\end{equation}
Also a straightforward computation shows that
\begin{equation}
  f''(\pi_1)= -\nu^2<0, \qquad f''\bigl(\frac1{\mu \pi_1}\bigr)
= \bigl( \gamma\nu \mu \pi_1(1-\pi_1) \bigr)^2 >0
\end{equation}

Due to the nature of the critical points, the point $z=\pi_1$ is
suitable for the steepest-descent analysis for $\mathcal{J}(v)$
and standard steepest-descent analysis will yield a good leading
term of the asymptotic expansion of $\mathcal{J}(v)$. However, for
$\mathcal{H}(u)$, the appropriate critical point is $z=1/(\mu
\pi_1)$, and in order to find the steepest-descent curve passing
the point $z=1/(\mu \pi_1)$, we need to deform the contour
$\Gamma$ through the pole $z=\pi_1$ and possibly some of
$\pi_{k+1}, \pi_{k+2}, \dots, \pi_{r}$. In the below, we will show
that the leading term of the asymptotic expansion of
$\mathcal{H}(u)$ comes from the pole $z=\pi_1$. Before we state
precise estimates, we first need some definitions.

Given any fixed $\epsilon >0$, we set
\begin{equation}\label{eq:qq2}
  \qq:= \pi_1 - \frac{\epsilon}{\nu \sqrt{M}}.
\end{equation}
Set
\begin{equation}
  \mathcal{H}_\infty(u):= ie^{-\epsilon u} \cdot \displaystyle \Res_{a=0}
\biggl( \frac1{a^k} e^{-\frac12 a^2-ua} \biggr), \qquad
  \mathcal{J}_\infty(v):= \frac{1}{2\pi} e^{\epsilon v}
\int_{\Sigma_\infty} s^k e^{\frac12s^2+vs} ds,
\end{equation}
where $\Sigma_\infty$ is the imaginary axis oriented from the
bottom to the top, and let
\begin{equation}
  Z_M:= \frac{(-1)^ke^{-Mf(\pi_1)}g(\pi_1)}{\nu^k M^{k/2}}.
\end{equation}

\begin{prop}\label{prop:Gauss}
Fix $\epsilon >0$ and set
$\qq$ by \eqref{eq:qq2}.
The followings hold for $M/N=\gamma^2$ with $\gamma$
in a compact subset of $[1,\infty)$.
\begin{itemize}
\item[(i)]
For any fixed $V\in\mathbb{R}$, there are constants $C, c>0$,
$M_0>0$ such that
\begin{equation}\label{eq:GaussJlimit}
  \biggl| \frac1{Z_M} \mathcal{J}(v) - \mathcal{J}_\infty(v) \bigr|
\le \frac{Ce^{-cv}}{\sqrt{M}}
\end{equation}
for $v\ge V$ when $M\ge M_0$.
\item[(ii)]
For any fixed $U\in\mathbb{R}$, there are constants $C, c>0$,
$M_0>0$ such that
\begin{equation}
  \bigl| Z_M \mathcal{H}(u) - \mathcal{H}_\infty(u) \bigr|
\le \frac{Ce^{-cu}}{\sqrt{M}}
\end{equation}
for $u\ge U$ when $M\ge M_0$.
\end{itemize}
\end{prop}

We prove this result in the following two subsections.

\subsection{Proof of Proposition \ref{prop:Gauss} (i)}\label{sec:Gauss1}

\begin{figure}[ht]
\centerline{\epsfxsize=11cm\epsfbox{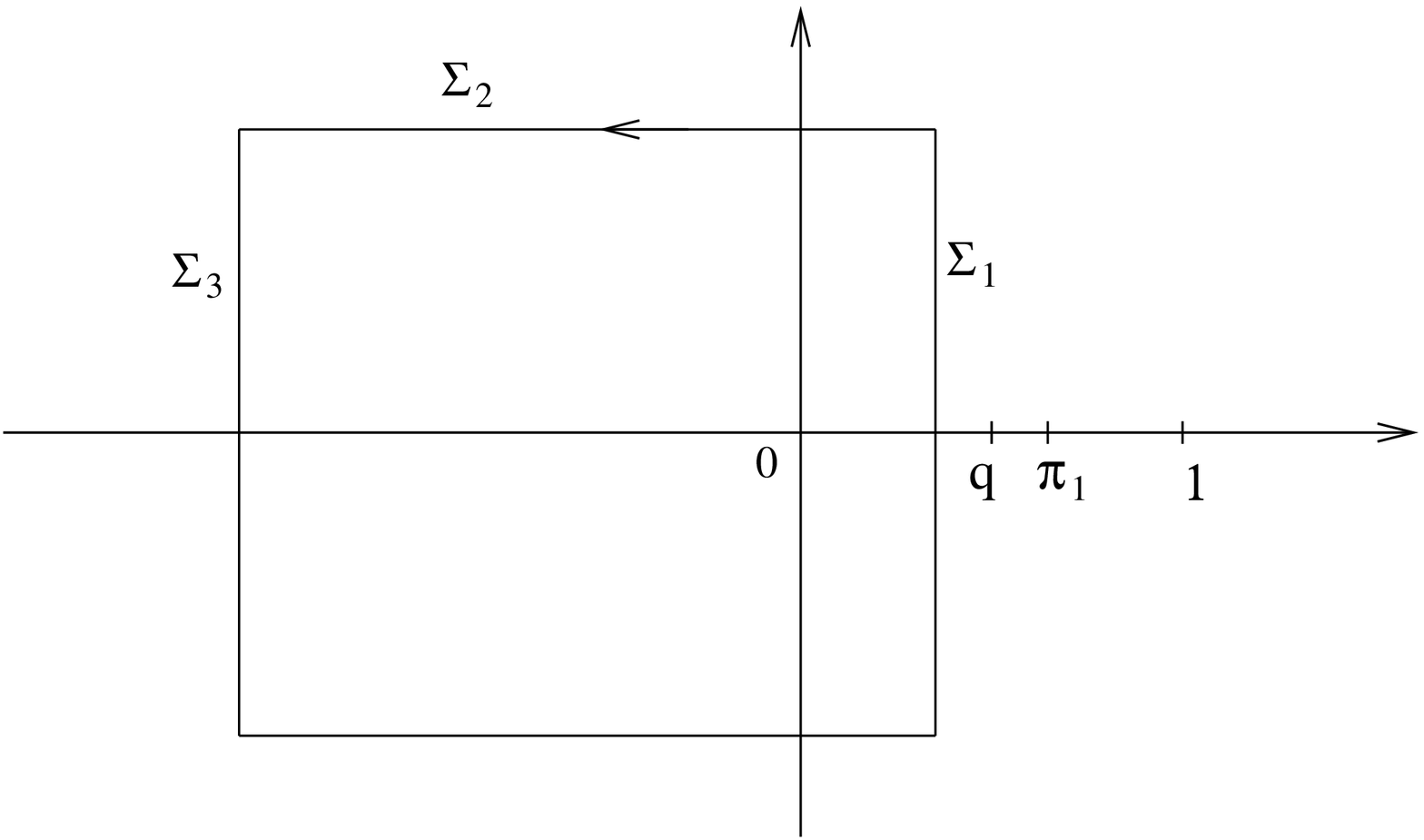}} \caption{Contour
$\Sigma$}\label{fig:Steep3}
\end{figure}

Let $R>0$ and define
\begin{eqnarray}
  \Sigma_1 &:=& \{ \pi_1-\frac{2\epsilon}{\nu\sqrt{M}}+iy : 0\le y\le 2  \} \\
  \Sigma_2 &:=& \{ \pi_1+2i -x : \frac{2\epsilon}{\nu\sqrt{M}}\le x\le R\} \\
  \Sigma_3 &:=& \{ \pi_1-R+i(2-y) : 0\le y\le 2 \},
\end{eqnarray}
and set
\begin{equation}
  \Sigma= \bigl( \cup_{k=1}^3 \Sigma_k \bigr)
\cup \overline{\bigl( \cup_{k=1}^3 \Sigma_k \bigr)}.
\end{equation}
The orientations of $\Sigma_j$, $j=1,2,3$ and $\Sigma$ are
indicated in Figure \ref{fig:Steep3}.

\begin{lemma}\label{lem:GaussJ}
For $\gamma\ge 1$, $Re(-f(z))$ is decreasing for $z\in \Sigma_1\cup\Sigma_2$
as $z$ travels on the contour along along
the prescribed orientation.
Also when $\gamma$ is in a compact subset of $[1,\infty)$, we can take $R>0$ large enough
so that
\begin{equation}\label{eq:bdS3}
  \max_{z\in\Sigma_3} Re(-f(z)) \le Re (-f(p_*)),
\end{equation}
where $p_*=\pi_1+2i$ is the intersection of $\Sigma_1$ and
$\Sigma_2$.
\end{lemma}

\begin{proof}
Any $z\in \Sigma_1$ is of the form $z=x_0+iy$, $0\le y\le 2$,
$x_0:=\pi_1-\frac{2\epsilon}{\nu\sqrt{M}}$. Set for $y\ge 0$,
\begin{equation}
  F_1(y):= Re(-f(x_0+iy))
= \mu(x_0-\qq) - \frac12 \ln(x_0^2+y^2) + \frac1{2\gamma^2}
\ln((1-x_0)^2+y^2).
\end{equation}
Then
\begin{equation}
  F_1'(y) = \frac{-y\bigl((\gamma^2-1)y^2+\gamma^2(1-x_0)^2-x_0^2\bigr)}
{\gamma^2(x_0^2+y^2)((1-x_0)^2+y^2)}.
\end{equation}
But as $0<x_0<\pi_1<\frac{\gamma}{1+\gamma}$, a straightforward
computation shows that $\gamma^2(1-x_0)^2-x_0^2
>0$. Therefore, $Re(-f(z))$
decreases as $z$ moves along $\Sigma_1$.

For $z\in\Sigma_2$, we have $z=\pi_1-x+2i$,
$\frac{2\epsilon}{\nu\sqrt{M}}\le x\le R$. Set
\begin{equation}
  F_2(x):= Re(-f(\pi_1-x+2i))
= \mu(\pi_1-\qq-x) - \frac12 \ln((\pi_1-x)^2+y^2) +
\frac1{2\gamma^2} \ln((1-\pi_1+x)^2+y^2).
\end{equation}
Then
\begin{equation}
  F_2'(x) = -\mu - \frac{x-\pi_1}{(x-\pi_1)^2+4} + \frac{x+1-\pi_1}{\gamma^2((x+1-\pi_1)^2+4)}.
\end{equation}
As the function $g(s)= \frac{s}{s^2+4}$ satisfies $-\frac14 \le g(s) \le \frac14$ for
all $s\in\mathbb{R}$, we find that for all $x\in\mathbb{R}$,
\begin{equation}
  F_2'(x) \le -\mu + \frac14 + \frac1{4\gamma^2}
= - \frac{4-\pi_1}{4\pi_1} - \frac{3+\pi_1}{4\gamma^2(1-\pi_1)}
\end{equation}
using the definition \eqref{eq:mu2} of $\mu$. But as
$0<\pi_1<\frac{\gamma}{\gamma+1}<1$, $F_2'(x)<0$ for all
$x\in\mathbb{R}$, and we find that $Re(-f(z))$ decreases as $z$
moves on $\Sigma_2$.

For $z\in\Sigma_3$, $z=\pi_1-R+i(2-y)$, $0\le y\le 2$. Then
for $\gamma$ in a compact subset of $[1,\infty)$, we can take $R>0$ sufficiently large
so that
\begin{equation}
\begin{split}
  F_3(y) &:= Re(-f(\pi_1-R+i(2-y)) \\
&= \mu (\pi_1-R-\qq) - \frac12 \ln((\pi_1-R)^2+(2-y)^2)
+ \frac1{2\gamma^2}\ln((1-\pi_1+R)^2+(2-y)^2)
\end{split}
\end{equation}
can be made arbitrarily small. However
\begin{equation}
  Re(-f(\pi_1+2i)) =
\mu (\pi_1-\qq) - \frac12 \ln(\pi_1^2+4)
+ \frac1{2\gamma^2}\ln((1-\pi_1)^2+4)
\end{equation}
is bounded for all $\gamma\ge 1$. Hence the result \eqref{eq:bdS3} follows.
\end{proof}

As $\gamma$ is in a compact subset of $[1,\infty)$, we assume that
\begin{equation}
  1\le \gamma\le \gamma_0
\end{equation}
for some fixed $\gamma_0\ge 1$. Also as $\pi_1$ is in a compact
subset of $(0,\frac{\gamma}{\gamma+1})$, we assume that there is
$0<\Pi<1/2$ such that
\begin{equation}
  \Pi\le \pi_1
\end{equation}
Fix $\delta$ such that
\begin{equation}
  0< \delta<\min \biggl\{ \frac{\Pi}2, \frac1{2(1+\gamma_0)^3}, \frac{\nu^2}{4C_1} \biggr\}, \qquad C_1:=
  \frac83\biggl( \frac1{\Pi^3}+(1+\gamma_0)^3\biggr) .
\end{equation}
Then for $|z-\pi_1|\le \delta$, by using the general inequality
\begin{equation}
\begin{split}
  \bigl| Re(-f(z)+f(\pi_1)+\frac12f''(\pi_1)(z-\pi_1) \bigr|
  &\le \biggl( \max_{|s-\pi_1|\le \delta} \frac1{3!}|f^{(3)}(s)|
  \biggr) |z-\pi_1|^3
\end{split}
\end{equation}
and the simple estimate for $|s-\pi_1|\le \delta$,
\begin{equation}
\begin{split}
  |f^{(3)}(s)| &= \biggl|
  \frac2{s^3}-\frac2{\gamma^2(s-1)^3}\biggr| \\
  &\le \frac2{(\pi_1-\delta)^3}+
  \frac2{\gamma_0^2(1-\pi_1-\delta)^3} \\
  &\le \frac{16}{\Pi^3}+ \frac{128}{\gamma_0^2} = 6C_1,
\end{split}
\end{equation}
we find that
\begin{equation}\label{eq:RefGaussJ}
\begin{split}
  \bigl| Re(-f(z)+f(\pi_1)+\frac12f''(\pi_1)(z-\pi_1) \bigr|
  &\le C_1 |z-\pi_1|^3 \\
  &\le \frac{\nu^2}{4} |z-\pi_1|^2, \qquad |z-\pi_1|\le \delta.
\end{split}
\end{equation}

 We split the contour
$\Sigma=\Sigma'\cup\Sigma''$ where $\Sigma'$ is the part of
$\Sigma$ in the disk $|z-\pi|\le \delta$, and $\Sigma''$ is the
rest of $\Sigma$. Let $\Sigma'_\infty$ be the image of $\Sigma'$
under the map $z\mapsto \nu\sqrt{M}(z-\pi_1)$ and let
$\Sigma''_\infty=\Sigma_\infty\setminus\Sigma_\infty'$. Set
\begin{equation}
  \mathcal{J}(v)=\mathcal{J}'(v)+\mathcal{J}''(v),
  \qquad \mathcal{J}_\infty(v)=\mathcal{J}_\infty'(v)+\mathcal{J}_\infty''(v)
\end{equation}
where $\mathcal{J}'(v)$ (resp. $\mathcal{J}_\infty'(v)$) is the
part of the integral formula of $\mathcal{J}(v)$ (resp.
$\mathcal{J}_\infty(v)$) integrated over the contour $\Sigma'$
(resp. $\Sigma'_\infty$).

Lemma \ref{lem:GaussJ} and the inequality \eqref{eq:RefGaussJ}
imply that
\begin{equation}
\begin{split}
  \max_{z\in\Sigma''} Re(-f(z)+f(\pi_1)) &\le Re(-f(z_0)+f(\pi_1))
  \\
  &\le Re(-\frac12f''(\pi_1)(z_0-\pi_1)^2)+
  \frac{\nu^2}4|z_0-\pi_1|^2 \\
  &= Re(\frac12\nu^2(z_0-\pi_1)^2)+
  \frac{\nu^2}4\delta^2.
\end{split}
\end{equation}
where $z_0$ is the intersection in the upper half plane of the
circle $|s-\pi_1|=\delta$ and the line
$Re(s)=\pi_1-\frac{2\epsilon}{\nu\sqrt{M}}$. As $M\to \infty$,
$z_0$ becomes close to $\pi_1+i\delta$. Therefore when $M$ is
sufficiently large,
\begin{equation}
  \max_{z\in\Sigma''} Re(-f(z)+f(\pi_1))
  \le -\frac{\nu^2}{12}\delta^2.
\end{equation}
Using this estimate and the fact that $Re(z-\pi_1)<0$ for
$z\in\Sigma''$, an argument similar to that in subsection
\ref{sec:Jmain} yields \eqref{eq:GaussJlimit}. We skip the detail.


\subsection{Proof of Proposition \ref{prop:Gauss} (ii)}\label{sec:Gauss2}

By using the Cauchy's residue theorem, for a contour $\Gamma'$
that encloses all the zeros of $g$ but $\pi_1$, we find
\begin{equation}
\begin{split}
  \mathcal{H}(u) =& i\nu \sqrt{M} \Res_{z=\pi_1} \biggl( e^{-\nu
  \sqrt{M}u(z-\qq)} e^{Mf(z)} \frac1{(\pi_1-z)^kg(z)}\biggr) \\
  &+ \frac{\nu \sqrt{M}}{2\pi} \int_{\Gamma'} e^{-\nu
  \sqrt{M}u(z-\qq)} e^{Mf(z)} \frac1{(\pi_1-z)^kg(z)} dz.
\end{split}
\end{equation}
Using the choice \eqref{eq:qq2} of $\qq$ and setting
$z=\pi_1+\frac{a}{\nu\sqrt{M}}$ for the residue term, we find that
\begin{equation}\label{eq:HGaussint}
\begin{split}
  Z_M \mathcal{H}(u) =& \mathcal{H}_1(u) + \frac{g(\pi_1)e^{-\epsilon u}}{2\pi (\nu\sqrt{M})^{k-1}} \int_{\Gamma'} e^{-\nu
  \sqrt{M}u(z-\pi_1)} e^{M(f(z)-f(\pi_1))} \frac{1}{(z-\pi_1)^kg(z)} dz.
\end{split}
\end{equation}
where
\begin{equation}
  \mathcal{H}_1(u):= ie^{-\epsilon u} \Res_{a=0} \biggl(
  \frac1{a^k}
  e^{-ua} e^{M\bigl( f(\pi_1+\frac{a}{\nu\sqrt{M}})-f(\pi_1)\bigr)}
  \frac{g(\pi_1)}{g(\pi_1+\frac{a}{\nu\sqrt{M}})}\biggr).
\end{equation}

We first show that $\mathcal{H}_1(u)$ is close to
$\mathcal{H}_\infty(u)$.
 Note that  all
the derivatives $f^{(\ell)}(\pi_1)$ and $g^{(\ell)}(\pi_1)$ are
bounded and $|g(\pi_1)|$ is strictly positive for $\gamma$ and
$\pi_1$ under our assumptions. The function
\begin{equation}
  e^{M\bigl( f(\pi_1+\frac{a}{\nu\sqrt{M}})-f(\pi_1)\bigr)}
\end{equation}
has the expansion of the form
\begin{equation}
  e^{-\frac12a^2 + a^2\bigl(c_1\bigl(\frac{a}{\sqrt{M}}\bigr)
  + c_2 \bigl(\frac{a}{\sqrt{M}}\bigr)^2+\cdots \bigr)}
\end{equation}
for some constants $c_j$'s when $a$ is close to $0$. On the other
hand, the function
\begin{equation}
   \frac{g(\pi_1)}{g(\pi_1+\frac{a}{\nu\sqrt{M}})}
\end{equation}
has the Taylor expansion of the form
\begin{equation}
  1+ c_1 \bigl(\frac{a}{\sqrt{M}}\bigr)
  + c_2 \bigl(\frac{a}{\sqrt{M}}\bigr)^2+\cdots
\end{equation}
for different constants $c_j$'s. Hence we find the expansion
\begin{equation}
\begin{split}
  &e^{-ua} e^{M\bigl( f(\pi_1+\frac{a}{\nu\sqrt{M}})-f(\pi_1)\bigr)}
  \frac{g(\pi_1)}{g(\pi_1+\frac{a}{\nu\sqrt{M}})} \\
  &\quad = e^{-ua -\frac12a^2} \biggl( 1+ \sum_{\ell, m=1}^\infty
  c_\ell a^{2\ell}\bigl(\frac{a}{\sqrt{M}}\bigr)^m
  + \sum_{\ell=1}^\infty d_\ell \bigl( \frac{a}{\sqrt{M}}\bigr)^\ell
  \biggr)
\end{split}
\end{equation}
for some constants $c_\ell, d_\ell$. Now as
\begin{equation}
  \Res_{a=0} \biggl( \frac1{a^\ell} e^{-au-\frac12a^2} \biggr)
\end{equation}
is a polynomial of degree at most $\ell-1$ in $u$, we find that
\begin{equation}
\begin{split}
  & \Res_{a=0} \biggl(
  \frac1{a^k}
  e^{-ua} e^{M\bigl( f(\pi_1+\frac{a}{\nu\sqrt{M}})-f(\pi_1)\bigr)}
  \frac{g(\pi_1)}{g(\pi_1+\frac{a}{\nu\sqrt{M}})}\biggr) \\
  &\quad =  \Res_{a=0} \biggl( \frac1{a^k}
  e^{-ua-\frac12 a^2}\biggr) + \sum_{j=1}^{k-1}
  \frac{q_j(u)}{(\sqrt{M})^j}
\end{split}
\end{equation}
for some polynomials $q_j$. Therefore, due to the factor
$e^{-\epsilon u}$ in $\mathcal{H}_1$, for any fixed
$U\in\mathbb{R}$, there are constants $C, c, M_0>0$ such that
\begin{equation}\label{eq:HGaussRes}
  \biggl| \mathcal{H}_1(u) - ie^{-\epsilon u} \Res_{a=0} \biggl(
  \frac1{a^k} e^{-ua-\frac12 a^2} \biggr) \biggr|
  \le \frac{Ce^{-cu}}{\sqrt{M}}
\end{equation}
for all $u\ge U$ when $M\ge M_0$.

\begin{figure}[ht]
\centerline{\epsfxsize=11cm\epsfbox{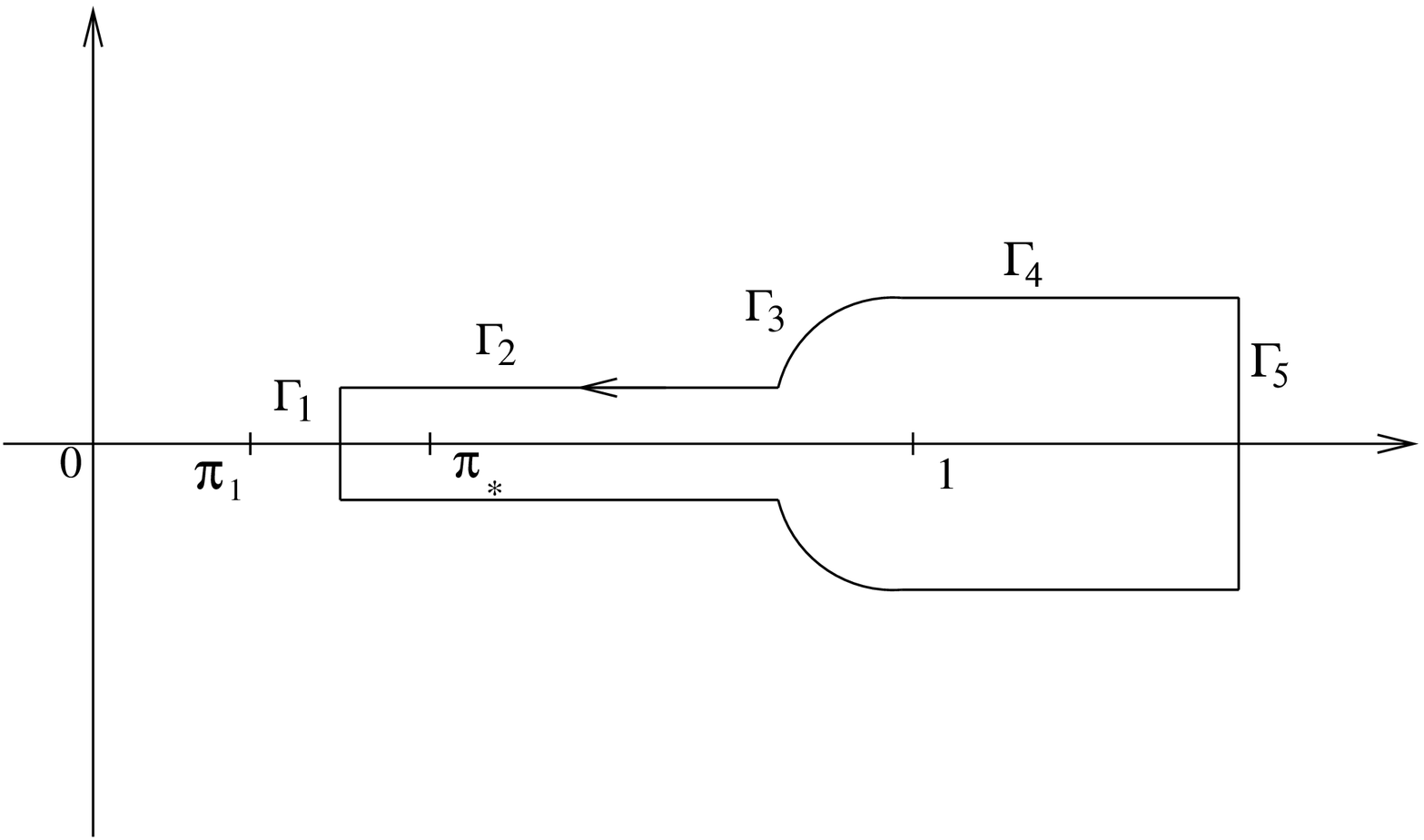}} \caption{Contour
$\Sigma'$}\label{fig:Steep4}
\end{figure}

Now we estimate the integral over $\Gamma'$ in
\eqref{eq:HGaussint}. We will choose $\Gamma'$ properly so that
the integral is exponentially small when $M\to\infty$. Let
$\pi_*:=\min\{ \pi_{k+1}\dots, \pi_r, 1, \frac1{\mu \pi_1}\}$.
Then $\pi_1=\dots=\pi_k<\pi_*$. Let $\delta>0$  and
$R>\max\{\pi_{k+1}, \dots, \pi_r, 1\}$ be determined in Lemma
\ref{lem:HGaussdecay} below. Define
\begin{eqnarray}
  \Gamma_1 &:=& \{ \frac{\pi_1+\pi_*}2+iy: 0\le y\le \delta \} \\
  \Gamma_2 &:=& \{ x+i\delta : \frac{\pi_1+\pi_*}2 \le x\le x_0  \} \\
  \Gamma_3 &:=& \{ 1+ \frac{1}{1+\gamma}e^{i(\pi-\theta)} : \theta_0\le \theta\le
  \frac{\pi}2 \} \\
  \Gamma_4 &:=& \{ x+i\frac1{1+\gamma} : 1 \le x\le R  \} \\
  \Gamma_5 &:=& \{ R+i(\frac1{1+\gamma}-y) : 0\le y\le \frac{1}{1+\gamma}  \}
\end{eqnarray}
where $x_0$ and $\theta_0$ are defined by the relation
\begin{equation}\label{eq:x0theta0}
  x_0+i\delta= 1+ \frac{1}{1+\gamma}e^{i(\pi-\theta_0)}.
\end{equation}
Set
\begin{equation}
  \Gamma'= \bigl(\cup_{j=1}^5 \Gamma_j \bigr) \cup \overline{\bigl(\cup_{j=1}^5 \Gamma_j
  \bigr)}.
\end{equation}
See Figure \ref{fig:Steep4} for $\Gamma'$ and its orientation.

\begin{lemma}\label{lem:HGaussdecay}
For $\gamma$ in a compact subset of $[1,\infty)$ and for $\pi_1$
in a compact subset of $(0,\frac{\gamma}{1+\gamma})$, there exist
$\delta>0$ and $c>0$ such that
\begin{equation}
  Re(f(z)-f(\pi_1))\le -c, \qquad z\in \Gamma_1\cup\Gamma_2.
\end{equation}
Also $Re(f(z))$ is a decreasing function in $z\in\Gamma_3
\cup\Gamma_4$, and when $R>\max\{\pi_{k+1},\dots, \pi_r,1\}$ is
sufficiently large,
\begin{equation}
  Re(f(z)-f(\pi_1)) \le Re\biggl(f(1+i\frac1{1+\gamma})-f(\pi_1)\biggr),
  \qquad z\in \Gamma_5.
\end{equation}
\end{lemma}

\begin{proof}
Note that
\begin{equation}
  |f'(z)| = \biggl|-\mu +\frac1{z}-\frac{\gamma^{-2}}{z-1}
  \biggr|
\end{equation}
is bounded for $z$ in the complex plane minus union of two compact
disks centered at $0$ and $1$. For $\gamma$ and $\pi_1$ under the
assumption, $\pi_1$ and $\frac{1}{\mu\pi_1}$ are uniformly away
from $0$ and $1$ (and also from $\frac{\gamma}{1+\gamma}$).
Therefore, in particular, there is a constant $C_1>0$ such that
for $z=x+iy$ such that $x\in [\pi_1,\frac1{\mu\pi_1}]$ and $y>0$,
\begin{equation}
  |Re(f(x+iy)-f(x))| \le |f(x+iy)-f(x)| \le \max_{0\le s\le 1} |f'(x+isy)|\cdot |y|
  \le C_1|y|.
\end{equation}
On the other hand, a straightforward calculation shows that when
$z=x$ is real, the function
\begin{equation}
  Re(f(x)-f(\pi_1)) = -\mu(x-\qq) + \ln|x| - \gamma^{-2}
  \ln|1-x| - f(\pi_1)
\end{equation}
decreases as $x$ increasing when $x\in (\pi_1, \frac{1}{\mu
\pi_1})$. (Recall that $\pi_1$ and $\frac1{\mu\pi_1}$ are the two
roots of $F'(z)=0$ and $0<\pi_1<\frac{\gamma}{1+\gamma} <
\frac1{\mu \pi_1}<1$.) Therefore we find that for $z=x+iy$ such
that $x\in [(\pi_1+\pi_{k+1})/2, \frac1{\mu\pi_1}]$,
\begin{equation}
\begin{split}
  Re(f(z)-f(\pi_1)) &\le Re(f(x)-f(\pi_1))+ C_1 |y| \\
  &\le Re \bigl(f(\frac{\pi_1+\pi_*}2) - f(\pi_1) \bigr)+ C_1 |y| .
\end{split}
\end{equation}
As $\pi_*$ is in a compact subset of $(\pi_1, 1)$, we find that
there is $c_1>0$ such that
\begin{equation}
  Re(f(z)-f(\pi_1)) \le -c_1 + C_1|y|
\end{equation}
for above $z$, and hence there are $\delta>0$ and $c>0$ such that
for $z=x+iy$ satisfying $|y|\le \delta$, $\frac{\pi_1+\pi_*}2\le
x\le \frac1{\mu \pi_1}$,
\begin{equation}\label{eq:ineqffc}
  Re(f(z)-f(\pi_1)) \le -c.
\end{equation}
Note that as $\frac1{\mu \pi_1}$ is in a compact subset of
$(\frac{\gamma}{\gamma+1}, 1)$ under out assumption, we can take
$\delta$ small enough such that $x_0$ defined by
\eqref{eq:x0theta0} is uniformly left to the point $\frac1{\mu
\pi_1}$. Therefore \eqref{eq:ineqffc} holds for
$z\in\Sigma_1\cup\Sigma_2$.

For $z=1+\frac1{1+\gamma}e^{i(\pi-\theta)}\in\Gamma_3$,
\begin{equation}
\begin{split}
  F_3(\theta) &:= Re\bigl(f(1+\frac1{1+\gamma}e^{i(\pi-\theta)})
  \bigr) \\
  &= -\mu(1+\frac1{1+\gamma}\cos(\pi-\theta) -\qq) + \frac12 \ln\bigl( 1+\frac2{1+\gamma} \cos(\pi-\theta)
  + \frac{1}{(1+\gamma)^2} \bigr) \\
  &\quad -\frac1{\gamma^2} \ln \biggl( \frac1{1+\gamma} \biggr).
\end{split}
\end{equation}
We set $t=\cos(\pi-\theta)$ and define
\begin{equation}
\begin{split}
  G(t)&:= F_3(\theta) =
  -\mu(1+\frac1{1+\gamma}t -\qq) + \frac12 \ln\bigl( 1+\frac2{1+\gamma} t
  + \frac{1}{(1+\gamma)^2} \bigr) +\frac1{\gamma^2} \ln(1+\gamma).
\end{split}
\end{equation}
Then
\begin{equation}
  G'(t) = -\mu\frac{1}{1+\gamma} + \frac{(1+\gamma)^{-1}}{1+2(1+\gamma)^{-1}t+(1+\gamma)^{-2}}
\end{equation}
is a decreasing function in $t\in[-1,1]$ and hence
\begin{equation}\label{eq:Gprimetet}
  G'(t)\le G'(-1)=\frac1{1+\gamma}\biggl( -\mu+\bigl(
  \frac{1+\gamma}{\gamma}\bigr)^2,
\biggr)\le 0
\end{equation}
as the function $\mu= \frac1{\pi_1}+\frac{\gamma^{-2}}{(1-\pi_1)}$
in $\pi_1\in (0, \frac{\gamma}{1+\gamma}]$ takes the minimum value
$\frac{(1+\gamma)^2}{\gamma^2}$ at $\pi=\frac{\gamma}{1+\gamma}$.
Therefore, $G(t)$ is a decreasing function in $t\in[-1,1]$ and
$Re(f(z))$ is a decreasing function in $z\in\Gamma_3$.

Set $y_1=\frac1{1+\gamma}$. For $z\in \Gamma_4$, $z=x+iy_1$, $x\ge
1$. Let
\begin{equation}
\begin{split}
  F_4(x) := Re \bigl( f(x+iy_1) \bigr)= -\mu (x-\qq) + \frac12 \ln (x^2+y_1^2) - \frac1{2\gamma^2} \ln
  \big( (x-1)^2+y_1^2 \bigr).
\end{split}
\end{equation}
Then
\begin{equation}
  F_4'(x) = -\mu + \frac{x}{x^2+y_1^2} - \frac{x-1}{\gamma^2\bigl( (x-1)^2+y_1^2
  \bigr)}.
\end{equation}
But the last term is non-negative and the middle term is less than
$1$ as $x\ge 1$. Also by the computation of \eqref{eq:Gprimetet},
$\mu \ge \frac{(1+\gamma)^2}{\gamma^2} \ge 1$. Therefore we find
that $F_4'(x) \le 0$ for $x\ge 0$, and $F_4(x)$ decreases as $x\ge
1$ increases.

Finally, for $z= R+iy$, $0\le y\le \frac{1}{1+\gamma}$,
\begin{equation}
  Re (f(R+iy) )
  = -\mu (R-\qq) + \frac12 \ln (x^2+y^2) - \frac1{2\gamma^2} \ln
  \big( (x-1)^2+y^2 \bigr)
\end{equation}
can be made arbitrarily small when $R$ is taken large enough,
while
\begin{equation}
  Re ( f(1+i\frac{1}{1+\gamma}) )
  = -\mu(1-\qq) + \frac12 \ln \biggl( 1+\frac1{(1+\gamma)^2}
  \biggr) - \frac1{\gamma^2} \ln \biggl( \frac1{1+\gamma} \biggr)
\end{equation}
is bounded.
\end{proof}

This lemma implies that
\begin{equation}
  Re(f(z)-f(\pi_1))\le -c
\end{equation}
for all $z\in\Gamma'$. Also note that $Re(z-\pi_1)>0$ for
$z\in\Gamma'$. Therefore for any fixed $U\in \mathbb{R}$, there
are constants $C, c, M_0>0$ such that
\begin{equation}
  \biggl| \frac{g(\pi_1)e^{-\epsilon u}}{2\pi (\nu\sqrt{M})^{k-1}}
\int_{\Gamma'} e^{-\nu
  \sqrt{M}u(z-\pi_1)} e^{M(f(z)-f(\pi_1))} \frac{1}{(z-\pi_1)^kg(z)} dz \biggr|
\le C e^{-\epsilon u} e^{-cM}
\end{equation}
for $M>M_0$ and for $u\ge U$. Together with \eqref{eq:HGaussRes},
this implies Proposition \ref{prop:Gauss} (ii).

\subsection{Proof of Theorem \ref{thm:main} (b)}

From the Proposition \ref{prop:Gauss} and the discussion in the
subsection \ref{sec:basicas}, we find that under the assumption of
Theorem \ref{thm:main} (b)
\begin{equation}\label{eq:proofmainblimit}
  \mathbb{P} \biggl( \bigr( \lambda_1 - \bigl(\frac1{\pi_1}+\frac{\gamma^{-2}}{1-\pi_1}\bigr) \bigr) \cdot
  \frac{\sqrt{M}}{\sqrt{\frac1{\pi_1^2}-\frac{\gamma^{-2}}{(1-\pi_1)^2}}} \le x
  \biggr)
\end{equation}
converges, as $M\to\infty$, to the Fredholm determinant of the
operator acting on $L^2((0,\infty))$ given by the kernel
\begin{equation}\label{eq:Proofmainbtemp}
  \int_0^\infty
  \mathcal{H}_\infty(x+u+y)\mathcal{J}_\infty(x+v+y)dy.
\end{equation}
Now we will express the terms $\mathcal{H}_\infty(u)$ and
$\mathcal{J}_\infty(v)$ in terms of the Hermite polynomials.

The generating function formula (see (1.13.10) of \cite{Koekoek})
of Hermite polynomials $H_n$,
\begin{equation}
  \sum_{n=0}^\infty \frac{H_n(x)}{n!} t^n = e^{2xt-t^2},
\end{equation}
implies that
\begin{equation}
  -ie^{\epsilon u} \mathcal{H}_\infty(u)=
  \displaystyle \Res_{a=0}
  \biggl( \frac1{a^k} e^{-\frac12 a^2-(u+y)a} \biggr)
  = \frac{(-1)^{k-1}(2\pi)^{1/4}}{\sqrt{(k-1)!}} p_{k-1}(u),
\end{equation}
where the orthonormal polynomial $p_{k-1}(x)$ is defined in
\eqref{eq:porthonomal}. The forward shift operator formula (see
(1.13.6) of \cite{Koekoek})
\begin{equation}
  H_n'(x)=2nH_{n-1}(x)
\end{equation}
implies that
\begin{equation}\label{eq:Gaussforward}
  p_k'(y)= \sqrt{k}p_{k-1}(y),
\end{equation}
and hence
\begin{equation}\label{eq:Hp}
  -ie^{\epsilon u} \mathcal{H}_\infty(u)=
  \displaystyle \Res_{a=0}
  \biggl( \frac1{a^k} e^{-\frac12 a^2-(u+y)a} \biggr)
  = \frac{(-1)^{k-1}(2\pi)^{1/4}}{\sqrt{k!}} p_k'(u).
\end{equation}
On the other hand, the Rodrigues-type formula (see (1.13.9) of
\cite{Koekoek})
\begin{equation}
  H_n(x) = e^{x^2} \biggl( -\frac{d}{dx} \biggr)^n [e^{-x^2}]
\end{equation}
implies that
\begin{equation}
  p_n(\xi)= \frac{(-1)^n}{(2\pi)^{1/4}\sqrt{n!}} e^{\xi^2/2}
  \biggl( \frac{d}{d\xi} \biggr)^n [e^{-\xi^2/2}].
\end{equation}
Sine the integral which appears in $\mathcal{J}_\infty(v)$ is
equal to
\begin{equation}
  \int_{\Sigma_\infty} s^k e^{\frac12s^2+vs} ds
  = \biggl(\frac{d}{dv}\biggr)^k \int_{\Sigma_\infty}
  e^{\frac12s^2+vs} ds
  = i \biggl(\frac{d}{dv}\biggr)^k e^{-v^2/2},
\end{equation}
we find
\begin{equation}\label{eq:Jp}
  e^{-\epsilon v} \mathcal{J}_\infty(v)
  = (-1)^k i (2\pi)^{-1/4} \sqrt{k!}  e^{-v^2/2} p_k(v).
\end{equation}

After a trivial translation, the Fredholm determinant of the
operator \eqref{eq:Proofmainbtemp} is equal to the Fredholm
determinant of the operator acting on $L^2((x,\infty))$ with the
kernel
\begin{equation}
\begin{split}
  K_2(u,v) &:=  \int_0^\infty
  \mathcal{H}_\infty(u+y)\mathcal{J}_\infty(v+y) dy.
\end{split}
\end{equation}
By \eqref{eq:Hp} and \eqref{eq:Jp},
\begin{equation}\label{eq:K2int}
\begin{split}
  (u-v) K_2(u,v) e^{\epsilon (u-v)}
  & = \int_0^\infty (u+y) \cdot p_k'(u+y) p_k(v+y) e^{-(v+y)^2/2} dy \\
  &\quad - \int_0^\infty p_k'(u+y) p_k(v+y) \cdot (v+y) e^{-(v+y)^2/2}
  dy.
\end{split}
\end{equation}
Note that $p_k$ satisfies the differential equation
\begin{equation}\label{eq:pode}
  p_k''(y)-yp_k'(y)+kp_k(y)=0,
\end{equation}
which follows from the differential equation
\begin{equation}
  H_n''(x) - 2x H_n'(x) + 2nH_n(x)=0
\end{equation}
for the Hermite polynomial (see (1.13.5) of \cite{Koekoek}). Now
use \eqref{eq:pode} for the first integral of \eqref{eq:K2int} by
and integrate by parts of the second integral by noting that
$(v+y)e^{-(v+y)^/2}= \frac{d}{dy} e^{-(v+y)^2/2}$ to obtain
\begin{equation}
\begin{split}
  (u-v)K_2(u,v) e^{\epsilon (u-v)}
  & = \int_0^\infty k p_k(u+y)p_k(v+y) e^{-(v+y)^2/2} dy \\
  & - p_k'(u)p_k(v)e^{-v^2/2}
  - \int_0^\infty p_k'(u+y)p_k'(v+y) e^{-(v+y)^2/2} dy .
\end{split}
\end{equation}
Note that the terms involving $p_k''(u+y)$ are cancelled out. Then
integrating by parts the second integral and noting that
$p_k'(u+y)=\frac{d}{dy} p_k(u+y)$, we obtain
\begin{equation}
  (u-v)K_2(u,v) e^{\epsilon (u-v)}  = - p_k'(u)p_k(v)e^{-v^2/2} + p_k(u)p_k'(v)
  e^{-v^2/2}.
\end{equation}
By using \eqref{eq:Gaussforward}, this implies that
\begin{equation}
  K_2(u,v) = \sqrt{k} e^{-\epsilon u}
  \frac{p_k(u)p_{k-1}(v)-p_{k-1}(u)p_k(v)}{u-v} e^{-v^2/2} e^{\epsilon
  v}.
\end{equation}
Therefore, upon conjugations, the Fredholm determinant of $K_2$ is
equal to the Fredholm determinant $\det(1-\mathbf{H}_x)$ in
\eqref{eq:Hop}. This completes the proof of Theorem \ref{thm:main}
(b).

\section{Samples of finitely many variables}\label{sec:finite}

In this section, we prove Proposition \ref{prop:finite}.

We take $M\to\infty$ and fix $N=k$. We suppose that
\begin{equation}
  \pi_1=\cdots=\pi_k.
\end{equation}
Then from \eqref{eq:eigendenoverU}, the density of the eigenvalues
is
\begin{equation}
  p(\lambda) = \frac1{C} V(\lambda)^2 \prod_{j=1}^k e^{-M\pi_1
  \lambda_j} \lambda_j^{M-k},
\end{equation}
and hence
\begin{equation}
  \mathbb{P}(\lambda_1\le t) = \frac1{C} \int_0^{t} \cdots \int_0^{t}
  V(y)^2 \prod_{j=1}^k e^{-M\pi_1
  y_j} y_j^{M-k} dy_j
\end{equation}
where 
\begin{equation}
  C= \frac{\prod_{j=0}^{k-1} (1+j)!(M-k+j)!}{(M\pi_1)^{Mk}}.
\end{equation}
By using the change of the
variables $y_j=\frac1{\pi_1}\bigl( 1+\frac{\xi_j}{\sqrt{M}}
\bigr)$,
\begin{equation}
  \mathbb{P} \biggl( \lambda_1\le \frac1{\pi_1}+ \frac{x}{\pi_1\sqrt{M}}
  \biggr)
  = \frac{e^{-kM}}{\pi_1^{kM}M^{k^2/2} C} \int_{-\sqrt{M}}^x\cdots \int_{-\sqrt{M}}^x V(\xi)^2
  \prod_{j=1}^k e^{-\sqrt{M}\xi_j}
  (1+\frac{\xi_j}{\sqrt{M}})^{M-k}
  d\xi_j .
\end{equation}
As $M\to\infty$ while $k$ is fixed, $\pi_1^{kM}M^{k^2/2}e^{kM} C
\to (2\pi)^{k/2}\prod_{j=0}^{k-1} (1+j)!$. By using the dominated
convergence theorem,
\begin{equation}\label{eq:finitetemp}
  \lim_{M\to\infty} \mathbb{P} \biggl( \bigl(\lambda_1-\frac1{\pi_1}
  \bigr) \pi_1\sqrt{M} x\biggr) = G_k(x).
\end{equation}

The result \eqref{eq:finiteinprob} follows from
\eqref{eq:finitetemp} and the fact that $G_k$ is a distribution
function.

\section{Last passage percolation, queueing theory
and heuristic arguments}\label{sec:last}

There is a curious connection between complex Gaussian sample
covariance matrices and a last passage percolation model.

Consider the lattice points $(i,j)\in \mathbb{Z}^2$. Suppose that
to each $(i,j)$, $i=1,\dots, N$, $j=1,\dots, M$, an independent
random variable $X(i,j)$ is associated. Let $(1,1) \nearrow (N,M)$
be the set of `up/right paths' $\pi=\{ (i_k, j_k)\}_{k=1}^{N+M-1}$
where $(i_{k+1},j_{k+1})-(i_k,j_k)$ is either $(1,0)$ or $(0,1)$,
and $(i_1,j_1)=(1,1)$ and $(i_{N+M-1}, j_{N+M-1})=(N,M)$. There
are $\binom{N+M-2}{N-1}$ such paths. Define
\begin{equation}\label{eq:lastppp}
  L(N,M) := \max_{\pi\in (1,1)\nearrow (N,M)} \sum_{(i,j) \in \pi}
  X(i,j).
\end{equation}
If $X(i,j)$ is interpreted as time spent to pass through the site
$(i,j)$, $L(N,M)$ is the \emph{last passage time} to travel from
$(1,1)$ to $(N,M)$ along an admissible up/right path.

Let $\pi_1, \dots, \pi_N$ be positive numbers. When $X(i,j)$ is the
exponential random variable of mean $\frac1{\pi_i M}$ (the density
function of $X(i,j)$ is $\pi_i M e^{-\pi_i M x}$, $x\ge 0$), it is
known that $L(N,M)$ has the same distribution as the largest
eigenvalue of the complex Gaussian sample covariance matrix of $M$
sample vectors of $N$ variables (see \eqref{eq:eigenden}): for
$M\ge N$,
\begin{equation}\label{eq:lastp}
  \mathbb{P}( L(N,M) \le x) = \frac1{C} \int_0^x \cdots \int_0^x
  \frac{\det\big( e^{-M\pi_i \xi_j} \bigr)_{1\le i,j\le
  N}}{V(\pi)} V(\xi) \prod_{j=1}^N \xi_j^{M-N} d\xi_j .
\end{equation}
We emphasize that $X(i,j)$, $j=1,2,\dots, M$ are identically distributed
for each fixed $i$. As a consequence, we have the following. Recall that
\begin{equation}
  \pi_j^{-1} = \ell_j.
\end{equation}

\begin{prop}
Let $L(M,N)$ be the last passage time in the above percolation model
with exponential random variables at each site. Let
$\lambda_1$ be the largest eigenvalue of $M$ (complex) samples of
$N\times 1$ vectors as in Introduction. Then for any $x\in \mathbb{R}$,
\begin{equation}
  \mathbb{P}(L(M,N)\le x) = \mathbb{P}(\lambda_1(M,N)\le x).
\end{equation}
\end{prop}

Formula \eqref{eq:lastp} for the case of $\pi_1=\cdots = \pi_N$
was obtained in Proposition 1.4 of \cite{kurtj:shape}. The general
case follows from a suitable generalization. Indeed, let $x_i,
y_j\in[0,1)$ satisfy $0\le x_iy_j<1$ for all $i,j$. When the
attached random variable, denoted by $Y(i,j)$, is the geometric
random variable of parameter $x_iy_j$ (i.e. $\mathbb{P}(X(i,j)=k)=
(1-x_iy_j)(x_iy_j)^k$, $k=0,1,2,\dots$), the last passage time,
$G(N,M)$, from $(1,1)$ to $(N,M)$ defined as in \eqref{eq:lastppp}
is known to satisfy
\begin{equation}\label{eq:Gss}
  \mathbb{P} ( G(N,M) \le n) = \prod_{i,j\ge 1}(1-x_iy_j) \cdot
  \sum_{\lambda : \lambda_1\le n}
  s_\lambda(x) s_\lambda(y)
\end{equation}
where the sum is over all partitions
$\lambda=(\lambda_1,\lambda_2,\dots)$ such that the first part
$\lambda_1\le n$, and $s_\lambda$ denotes the Schur function,
$x=(x_1,x_2, \dots)$ and $y=(y_1,y_2,\dots)$. This identity was
obtained by using the Robinson-Schensted-Knuth correspondence
between generalized permutations (matrices of non-negative integer
entries) and pairs of semistandard Young tableaux (see, e.g.
\cite{kurtj:shape}, \cite{Okounkov}, (7.30) of \cite{BR1}). The
normalization constant follows from the well-known Cauchy identity
(see, e.g. \cite{Stl})
\begin{equation}
  \sum_{\lambda} s_\lambda(x) s_\lambda(y) = \prod_{i,j\ge 1}
  (1-x_iy_j)
\end{equation}
where the sum is over all partitions.  Now set
$x_i=1-\frac{M\pi_i}{L}$, $i=1,\dots, N$, $x_i=0, i>N$ and
$y_j=1$, $j=1,\dots, M$, $y_j=0, j>M$, and $n=xL$. By taking
$L\to\infty$, it is easy to compute that $\frac1{L} Y(i,j)$
converges to the exponential random variable $X(i,j)$, while one
can check that the summation on the right-hand-side of
\eqref{eq:Gss} converges to the right-hand-side of
\eqref{eq:lastp}, and hence the identity \eqref{eq:lastp} follows.
There are determinantal formulas for the right-hand-side of
\eqref{eq:Gss} (see e.g. \cite{Gessel}, \cite{kurtj:shape},
\cite{Okounkov}), some of which, by taking the above limit, would
yield an alternative derivation of Proposition
\ref{prop:integralformula}.

\bigskip

Another equivalent model is a queuing model. Suppose that there
are $N$ tellers and $M$ customers. Suppose that initially all $M$
customers are on the first teller in a queue. The first customer
will be served from the first teller and then go to the second
teller. Then the second customer will come forward to the first
teller. If the second customer finishes his/her business before
the first customer finishes his/her business from the second
teller, the second customer will line up a queue in the second
teller, and so on. At any instance, only one customer can be
served at a teller and all customers should be served from all
tellers in the order. The question is the total exit time $E(M,N)$
for $M$ customers to exit from $N$ queues. We assume that the
service time at teller $i$ is given by the exponential random
variable of mean $\frac1{\pi_i M}$. Assuming the independence,
consideration of the last customer in the last queue will yield
the recurrence relation
\begin{equation}
  E(M,N) = \max\{ E(M-1, N), E(M,N-1)\} + e(N)
\end{equation}
where $e(N)$ denotes the service time at the teller $N$. But note that
the last passage time in the percolation model also satisfies the same recurrence
relation
\begin{equation}
  L(M,N) = \max\{ L(M-1, N), L(M,N-1)\} + X(M,N)
\end{equation}
where $X(M,N)$ is the same exponential random variable as $e(N)$.
Therefore we find that $E(M,N)$ and $L(M,N)$ have the same distribution.
Thus all the results in Introduction also applied to $E(M,N)$.

\bigskip

Now we indicate how the critical value $\ell_j= \pi_j^{-1}= 1+\gamma^{-1}$ of
Theorem \ref{thm:main} can be predicted in the last passage
percolation model.

First, when $X(i,j)$ are all identical exponential random
variables of mean $1$, Theorem 1.6 of \cite{kurtj:shape} shows
that, as $M,N\to\infty$ such that $M/N$ is in a compact subset of
$(0,\infty)$, $L(N,M)$ is approximately
\begin{equation}\label{eq:lp1}
  L(N,M) \sim L(M,N) \sim (\sqrt{M}+\sqrt{N})^2 +
  \frac{(\sqrt{M}+\sqrt{N})^{4/3}}{(MN)^{1/6}}  \chi_0
\end{equation}
where $\chi_0$ denotes the random variable of the GUE Tracy-Widom
distribution. On the other hand, note that when $N=1$, $X(1,j)$
are independent, identically distributed exponential random
variables of mean $\frac1{M\pi_1}$, and hence the classical
central limit theorem implies that the last passage time from
$(1,1)$ to $(1,xM)$ is approximately
\begin{equation}\label{eq:lp2}
  \pi_1^{-1}x + \frac{x\pi_1^{-1}}{\sqrt{M}} g
\end{equation}
where $g$ denotes the standard normal random variable. Note the
different fluctuations which is due to different dimensions of two
models.

\begin{figure}[ht]
\centerline{\epsfxsize=8cm\epsfbox{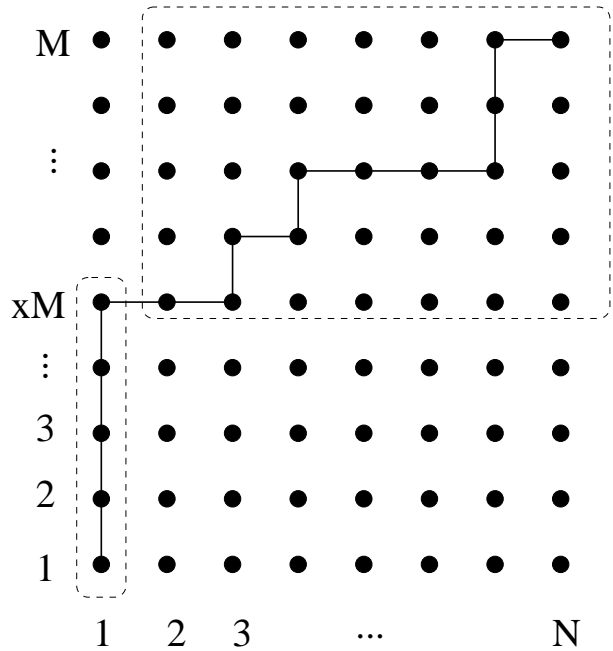}} \caption{Last
passage percolation when $r=1$} \label{fig:last}
\end{figure}

Now consider the case when $r=1$ in Theorem \ref{thm:main} i.e.
$\pi_2=\pi_3=\dots = \pi_N=1$; $X(i,j)$, $1\le j\le M$ is
exponential of mean $\frac1{M\pi_1}$ and $X(i,j)$, $2\le i\le N$,
$1\le j\le M$ is exponential of mean $\frac1{M}$. We take
$M/N=\gamma^2\ge 1$. An up/right path consists of two pieces; a
piece on the first column $(1,j)$ and the other piece in the
`bulk', $(i,j)$, $i\ge 2$. Of course the first part might be
empty. We will estimate how long the last passage path stays in
the first column. Consider the last passage path
\emph{conditioned} that it lies on the first column at the sites
$(1,1), (1,2),\dots (1,xM)$ and then enters to the bulk $(i,j),
i\ge 2$. See Figure \ref{fig:last}.
Then from \eqref{eq:lp2} and \eqref{eq:lp1}, we expect that the
(conditioned) last passage time is, to the leading order,
\begin{equation}
  f(x) = \pi_1^{-1}x + (\sqrt{1-x} +\gamma^{-1})^2.
\end{equation}
It is reasonable to expect that the last passage time is the maximum of $f(x)$
over $x\in [0,1]$, to the leading order. An elementary Calculus
shows that
\begin{equation}\label{eq:lpf}
  \max_{x\in [0,1]} f(x) = \begin{cases}
  f(0) = (1+\gamma^{-1})^2, \qquad &\text{if $\pi^{-1}\le
  1+\gamma^{-1}$} \\
  f\bigl( 1- \frac{\gamma^{-2}}{(\pi^{-1}-1)^2} \bigr)
  = \frac1{\pi_1}+\frac{\gamma^{-2}}{1-\pi_1},
  \qquad &  \text{if $\pi^{-1} >  1+\gamma^{-1}$.}
  \end{cases}
\end{equation}
When $\max f$ occurs at $x=0$, the last passage path enters
directly into the bulk and hence the fluctuation of the last
passage time is of of $M^{-2/3}$ due to \eqref{eq:lp1}. But if the
$\max f(x)$ occurs for some $x>0$, then the fluctuation is
$M^{-1/2}$ due to \eqref{eq:lp2}, which is larger than the
fluctuation $M^{-2/3}$ from the bulk. Note that the value of $\max
f$ in \eqref{eq:lpf} agrees with the leading term in the scaling
of Theorem \ref{thm:main}. This provides an informally explanation
of the critical value $1+\gamma^{-1}$ of $\pi_1^{-1}$.

If one can make this kind of argument for the sample covariance
matrix, one might be able to generalize it to real sample
covariance matrix.


\end{document}